\documentclass[11pt,a4paper,reqno]{amsart}
\usepackage{import}
\usepackage{p00preamble}
\usepackage[
    citestyle       =   numeric-comp,
    style           =   ext-numeric,
    date            =   year,
    articlein       =   false,
    giveninits      =   true,
    abbreviate      =   true,
    doi             =   false,
    isbn            =   false,
    url             =   false,
    maxbibnames     =   6,
    maxcitenames    =   6,
    sorting         =   nyt,
]{biblatex}
\usepackage{p01bibformat}
\usepackage{geometry}
\title{Trace type Orlicz spaces and analysis of Orlicz spaces by Lebesgue exponents}
\author{Albin Petersson }
\date{July 4, 2025}
\address{Department of Mathematics,
Linn{\ae}us University, V{\"a}xj{\"o}, Sweden}
\email{albin.petersson@lnu.se}

\begin{document}
\begin{abstract}
In the paper, we analyze the Lebesgue exponents $p_\Phi$ and $q_\Phi$, and show that for any $p_\Phi< p < \infty$ and $1< q<q_\Phi$, there exists an equivalent Young function $\Psi$ with $p < p_\Psi < \infty$ and $1<q_\Psi < q$. This type of construction is used to improve upon the inclusions $L^{p_\Phi}\cap L^{q_\Phi}\subseteq L^\Phi \subseteq L^{p_\Phi} + L^{q_\Phi}$. For trace type Orlicz spaces $L^{\Phi,\Phi}$, we find that when $\Phi \in \Delta_2$, we have $L^{\Phi,\Phi} \subseteq L^\Phi$ if and only if $\Phi(\nm{f}{L^\Phi}) \le C \rho_\Phi(f)$ for all $f\in L^\Phi$, and the reverse inclusion is equivalent to the reversed inequality.
\end{abstract}
\maketitle
% ------------- Title Page -------------
%\includepdf[pages={1}]{coversheet.pdf}
\thispagestyle{empty}

% --------- Table of Contents ---------
\newpage
\thispagestyle{empty}
\tableofcontents
% ----------------------------------------
%  Introduction              
% ----------------------------------------
\newpage
\setcounter{page}{1}
\section{Introduction}\label{sec:intro}
A fundamental setting for mathematics in numerous fields is 
the Lebesgue space $L^p$, $1\le p \le \infty$ 
(with important special cases $p=1, 2$, or $\infty$). 
However, some problems require the examination of 
non-homogeneous functionals which do not neatly slot into any 
Lebesgue space. One example from statistical physics is 
the entropy $E$ of a probability density function $f$ given by
$$
E(f) = - \int f(\xi) \log f(\xi) \, d\xi,
$$
which is not well-defined for all $f\in L^1$. 
To solve this problem, one turns to the spaces 
introduced by W. Orlicz in \autocite{orli:1932}, which 
we now call Orlicz spaces. These are denoted by $L^\Phi$, 
where $\Phi$ is a Young function, a certain type of convex function. 
The special case $\Phi(t) = t^p$, $1\le p < \infty$, gives 
$L^\Phi = L^p$. By \autocite{maje:2013}, one can avoid 
the continuity problem for $E$ above by replacing 
the classical statistical framework $\scal{L^\infty}{L^1}$ with $\scal{L^\Psi}{L^\Phi}$, where 
$\Phi(t) = t\log(1+t)$ and $\Psi(t) = \cosh t - 1$.

Since the Orlicz space $L^\Phi$ is completely determined by 
the Young function $\Phi$, it is clear that 
the analysis of Young functions and their properties is essential. 
This analysis is mainly done via certain indices, 
which provide an estimate of the growth rate of 
a Young function $\Phi$. We consider quantities such as 
the limits at the origin and infinity of 
$\frac{t \Phi'(t)}{\Phi(t)}$ and $\frac{\ln \Phi(t)}{\ln t}$. 
Such indices, and many more, 
appear in \autocite{mali:1985,mali:1989B}.

Much of the analysis of this paper revolves around 
the so-called Lebesgue exponents $p_\Phi$ and $q_\Phi$, 
which appeared in \autocite{simo:1964}, and have subsequently been analyzed in, e.g., \autocite{boni:2024}. We will mainly focus on 
the case that $p_\Phi < \infty$, 
which is equivalent to 
$\Phi$ fulfilling the (global) $\Delta_2$-condition 
(see \cref{rem:delta2}).

In this paper, we examine the Lebesgue exponents closely and 
determine their most glaring limitations. 
Namely, we show how equivalent Young functions can be constructed 
with widely differing Lebesgue exponents. 
We also present some alternative methods of obtaining 
growth estimates of Young functions, and by extension, 
obtain inclusion properties for Orlicz spaces. 
Furthermore, we investigate mixed Orlicz spaces, 
which are useful in interpolation theory. 
Particularly, we investigate the 
trace type Orlicz spaces $L^{\Phi,\Phi}(\rdd)$, 
where the same Young function $\Phi$ is used for 
the two $d$-dimensional variables, but separately, 
as opposed to $L^\Phi(\rdd).$

In the case of Orlicz modulation spaces, 
as was shown in \autocite[Corollary 2.5]{toft:2023}, 
we have $M^\Phi(\rd) = M^{\Phi,\Phi}(\rd)$, 
meaning that $\nm{V_\fy \cdot}{L^{\Phi}}$ 
is equivalent to 
$\nm{V_\fy \cdot}{L^{\Phi,\Phi}}$, 
where $V_\fy$ denotes the short-time Fourier transform 
with window function $\fy$. 
This equivalence does not hold in general.

Inclusions between the spaces $L^{\Phi,\Phi}$ and $L^\Phi$ 
are shown to be linked to inequalities between 
the quantities 
$$
\rho_\Phi(f) =\int \Phi(|f(x)|)\,dx
$$ 
and $\Phi(\nm{f}{L^\Phi})$. 
We find that for $\Phi\in \Delta_2$, 
$L^{\Phi,\Phi} \subseteq L^\Phi$ 
if and only if 
$\rho_\Phi(f)\le C \Phi(\nm{f}{L^\Phi})$ 
for some constant $C>0$, which in turn is equivalent to 
$\Phi$ being submultiplicative. The reversed inclusion is equivalent to the reverse inequality holding, which in turn is equivalent to  $\Phi$ being supermultiplicative. Evidently, $L^\Phi = L^{\Phi,\Phi}$ therefore holds if and only if 
$$
C^{-1} \Phi(\nm{f}{L^\Phi})\le\rho_\Phi(f)\le C \Phi(\nm{f}{L^\Phi})
$$
for some constant $C>0$, which by \autocite{fino:1991} is also equivalent to $L^\Phi = L^p$ for some $p\ge 1$.

The paper is organized as follows.
In \cref{sec:prelim}, we introduce necessary notations and definitions. The aforementioned Lebesgue exponents are 
examined in \cref{subsec:lebexp}. 
Their limitations are described in \cref{subsec:expo}, 
where it is found that if $p_\Phi < p < \infty$ and $1<q < q_\Phi$, 
then there exists an equivalent Young function $\Psi$ 
with $p<p_\Psi < \infty$ and $1<q_\Psi < q$ (\cref{prop:p1p2LPhi}). 
In \cref{subsec:equivYfn} and \cref{subsec:different}, 
we rectify these problems in different ways by focusing on shared properties for the equivalence class and ultimately obtain better growth estimates for $\Phi$ (\cref{prop:inclusionsEquiv} and \cref{thm:lnPhiinclusions}).

The relation between spaces $L^\Phi$ and $L^{\Phi,\Phi}$ 
is explored in \cref{sec:phi}. Here, \cref{subsec:modular} 
is focused on the Young function modular 
$\rho_\Phi(f) = \int \Phi(|f(x)|)\, dx$, and 
\cref{subsec:submult} contains sufficient and necessary conditions 
on $\Phi$, $\Phi_1$, and $\Phi_2$ for inclusions between $L^\Phi$ and $L^{\Phi_1,\Phi_2}$ to hold. In particular, such results are stated in the case that $\Phi=\Phi_1=\Phi_2$ (\cref{cor:tracetype} and \cref{cor:submult}).
% ----------------------------------------
%  Preliminaries            
% ----------------------------------------
\section{Preliminaries}\label{sec:prelim}
By $A\subseteq B$, we mean that $A$ is a (not necessarily proper) subset of $B$, and we will write $A\subsetneq B$ to specify that $A$ is a proper subset of $B$. We will use $d$ exclusively to denote dimension, and we will deal only with functions whose domains are 
$$
\re_+ = \sets{t\in\re}{t\ge0},
$$ 
$\rd$, or $\rdd$.

For $p\in[1,\infty]$, let $L^p(\rd)= L^p$ denote the usual Lebesgue space with norm
$$
\nm{f}{L^p(\rd)} \equiv \nm{f}{L^p} \equiv \nm{f}{p} \equiv
\begin{cases} 
\left( \int_{\rd} |f(x)|^p, d x \right)^{\frac 1 p}, & p<\infty, 
\\[1ex]
\displaystyle{\esssup_{x\in \rd}} |f(x)|, & p=\infty, 
\end{cases}
$$
where $\defrd$ is a Lebesgue measurable function. The norm $\nm{\cdot}{L^p}$ imposes certain decay conditions on the functions in $L^p$ which depend on the parameter $p$. Using different such conditions with respect to different variables, we arrive at the definition of mixed norm Lebesgue spaces, which we recall below.

\begin{defn}
The \emph{mixed Lebesgue space} $L^{p,q}\left(\rdd\right)$ consists of all Lebesgue measurable functions $\defrdd$ such that
$$
 \nm{f}{L^{p,q}(\rdd)} \equiv\nm{f}{L^{p,q}} \equiv\nm{f}{p,q}
\equiv 
\nm{f_2}{q}
$$
is finite, where
$$
f_2(y) = \nm{f(\cdot,y)}{p}.
$$
\end{defn}

Note that if $p=q$ in the definition above, then $L^{p,p}(\rdd) = L^p(\rdd)$.

We recall the following definition of Young function.
\begin{defn}\label{Def:YoungFunc}
A function $\fdef{\Phi}{[0,\infty)}{[0,\infty]}$ is called a \emph{Young function} if
\begin{enumerate}
  \item \label{def:YFitem1}$\Phi$ is convex,

  \vrum

  \item \label{def:YFitem2}$\Phi (0)=0$,

  \vrum

  \item \label{def:YFitem3}$\Phi (t)<\infty$ for some $t>0$,

  \vrum

  \item \label{def:YFitem4}$\lim
\limits_{t\rightarrow\infty} \Phi (t)=+\infty$.
\end{enumerate}
\end{defn}

A Young function does not have to be continuous by this definition. In fact, we do not exclude $\Phi$ with infinite function values, hence the function
$$
\Phi(t) = 
\begin{cases}
0 & 0\le t<1, \\
t-1 & 1\le t < 2, \\
\infty & t \ge 2
\end{cases}
$$
is a Young function. However, the only possible discontinuities for a Young function $\Psi$ appear at points $t_0>0$ such that $\Psi(t) =\infty$ for all $t>t_0$, since otherwise $\Psi$ would not be convex. Evidently, this means that $\Psi$ can have at most one jump.

The reason we require \ref{def:YFitem3} in the definition above is so that each Young function has a well-defined complementary Young function, given by
$$
\Psi(t) = \sup_{s\ge0} \left(s t - \Phi(s) \right).
$$

Using the notations of \autocite[Section 1]{boni:2024}, we recall the definition of Lebesgue exponent as follows.

\begin{defn}
Let $\Phi$ be a Young function and let $\Omega = \sets{t>0}{0<\Phi(t)<\infty}$. Then the \emph{Lebesgue exponents} $p_\Phi$ and $q_\Phi$ are given by
\begin{align}
p_{\Phi}
&\equiv
\begin{cases}
{\displaystyle{\sup_{t\in \Omega}
\left (\frac{t \Phi_+ '(t)}{\Phi(t)}\right )}},
& \Omega=(0,\infty),
\\[1ex]
\infty, & \Omega\neq (0,\infty),   
\end{cases}
\label{Eq:StrictYoungFunc2}
\intertext{and}
q_{\Phi}
&\equiv
\begin{cases}
{\displaystyle{\inf_{t\in \Omega}
\left (\frac{t \Phi_+'(t)}{\Phi(t)}\right )}},
& \Omega\neq \emptyset,
\\
\infty, & \Omega = \emptyset.
\end{cases}
\label{Eq:StrictYoungFunc1}
\end{align}
\end{defn}

Since, for any $a>0$, $\Phi'(a)$ is the slope of the tangent to $\Phi(t)$ in the point $t=a$ (or the slope of the right tangent if $\Phi$ is not differentiable in $a$), and $\Phi(a)/a$ is the slope of the line segment between the origin and $(a,\Phi(a))$, it is clear that $a\Phi'(a)/\Phi(a)\ge 1$ for all $a>0$ by the convexity of $\Phi$. Hence, for all Lebesgue exponents we must have $1\le q_\Phi \le p_\Phi \le \infty$. An important special case is when $q_\Phi = p_\Phi = p$, as then $\Phi(t) = C t^p$ for some constant $C>0$, and $L^\Phi = L^p$.

The Lebesgue exponents may be found in \autocite{simo:1964}, which is their first appearance according to \autocite{mali:1985}. These indices and similar ones of note for the study of Orlicz spaces can be found in \autocite[p. 20]{mali:1985}.

We will often require that a Young function $\Phi$ fulfills $p_\Phi < \infty$. Since this is equivalent to $\Phi$ fulfilling the $\Delta_2$-condition (see \cref{rem:delta2}), we write $\Phi\in \Delta_2$ whenever $\Phi$ is a Young function with $p_\Phi < \infty$. If $\Phi\in\Delta_2$ then $\Phi$ fulfills $0<\Phi(t)<\infty$ for all $t>0$, but this condition is not sufficient to guarantee that $p_\Phi<\infty$. For example, if $\Phi(t) = e^t - 1$, then $\Phi$ is a Young function with $0<\Phi(t)<\infty$ for all $t>0$, but $p_\Phi = \infty$.  We will also consider the \qu{nice} Young functions of \autocite{rao:1991}, the so-called \emph{N functions}. Recall that a Young function $\Phi$ is an N function if
$$
\lim_{t\rightarrow 0^+} \frac{\Phi(t)}{t} = 0 \quad\text{and}\quad
\lim_{t\rightarrow \infty} \frac{\Phi(t)}{t} = \infty.
$$

\begin{rem}
Young functions are not necessarily differentiable, but being convex, they are still semi-differentiable. Since the definition above, as well as all computations and results involving derivatives of Young functions in this paper, are the same whether one uses the left derivative or the right derivative, we will simply perform these calculations with the right derivative, and let $\Phi'\equiv \Phi'_+$ unless otherwise stated. 
\end{rem}

\begin{rem}\label{rem:qphi}
As the Lebesgue exponents will be central to our discussions in Section \ref{subsec:lebexp}, we will mention some basic properties, (cf. \autocite{boni:2024}). Observe that if $p_\Phi < \infty$, there are constants $C_1,C_2>0$ such that 
$$
C_1 t^{p_\Phi} \le \Phi(t) \le C_2 t^{q_\Phi}
$$
for $t\leq 1$, and
$$
C_1 t^{q_\Phi} \le \Phi(t) \le C_2 t^{p_\Phi}
$$
for $t\ge 1$. It follows that if $q_\Phi = p_\Phi$, then
$$
\Phi(t) = C t^{p_\Phi}
$$
for some constant $C>0$. Moreover, if
$$
\Phi_1(t) = \begin{cases}
C_1t^{p_\Phi} & t\le 1 \\
C_1 t^{q_\Phi} & t > 1,
\end{cases}
\quad\text{and}\quad \Phi_2(t) = \begin{cases}
C_2t^{q_\Phi} & t\le 1 \\
C_2 t^{p_\Phi} & t > 1,
\end{cases}
$$
then $L^{\Phi_1} = L^{p_\Phi} + L^{q_\Phi}$, $L^{\Phi_2} = L^{p_\Phi} \cap L^{q_\Phi}$,
and there exist constants $C_3, C_4>0$ such that
$$
\Phi_1(t) \le C_3 \Phi(t) \le C_4 \Phi_2(t), \quad t\ge 0,
$$
hence
\begin{equation}\label{eq:qPhiincl}
L^{p_\Phi} \cap L^{q_\Phi} \subseteq L^{\Phi} \subseteq L^{p_\Phi} + L^{q_\Phi}.    
\end{equation}
\end{rem}

To simplify notations,  we will let
$\rho_\Phi(f) = \nm{\Phi(|f|)}{L^1}$, and we will refer to $\rho_\Phi$ as \emph{the $\Phi$ modular}.

\begin{defn}
Let $\Phi $ be a Young function.
The \emph{Orlicz space} $L^{\Phi }(\rd)$ consists
of all Lebesgue measurable functions
$\defrd$ such that
$$
\nm{f}{L^{\Phi}(\rd)}
\equiv
\nm{f}{L^\Phi}
\equiv
\nm{f}{\Phi}
\equiv
\inf 
\Sets{\lambda>0}{\rho_\Phi\Big(\frac{f}{\lambda}\Big)\leqslant 1}
$$
is finite.
\end{defn}

In \cref{subsec:lebexp}, we will need to discuss Young functions which are equivalent in the sense that they generate the same Orlicz space. For this reason, we recall the following definition.

\begin{defn}
Let $f$ and $g$ be non-negative functions defined on $\re_+$. Then $f$ and $g$ are said the be \emph{equivalent} if there exists a constant $C>0$ such that
$$
C^{-1} g(t) \le f(t) \le Cg(t), \quad t\ge 0.
$$
\end{defn}

In particular, this defines an equivalence relation that partitions the set of Young functions into equivalence classes.
Note that if two Young functions $\Phi$ and $\Psi$ are equivalent, then so too are their norms, and therefore $L^{\Phi} = L^{\Psi}$. Hence an equivalence class of Young functions, which we will denote $[\Phi]$, is the set of all Young functions $\Psi$ for which $L^\Psi = L^\Phi$.

Of particular note in \cref{sec:phi} are the mixed Orlicz spaces, whose definition we recall as follows.

\begin{defn}
Let $\Phi,\Psi$ be Young functions. The \emph{mixed Orlicz space} $L^{\Phi,\Psi}(\rdd)$ consists of all Lebesgue measurable functions $\defrdd$ such that
$$
\nm{f}{L^{\Phi,\Psi}(\rdd)} \equiv \nm{f}{L^{\Phi,\Psi}} \equiv \nm{f}{\Phi,\Psi} \equiv \Nm{f_2}{L^{\Psi}(\rd)}
$$
is finite, where
$$
f_2(y) = \nm{f(\cdot,y)}{L^{\Phi}(\rd)}.
$$
\end{defn}

In this paper, we will focus on the case that $\Phi = \Psi$ in the definition above, and refer to these as \emph{trace type Orlicz spaces}. For $\Phi(t) =t^p$, it is clear that $L^{\Phi,\Phi}(\rdd) = L^\Phi(\rdd)$.

Lastly, we will recall the following definition.
\begin{defn}
A function $\fdef{\Phi}{[0,\infty)}{[0,\infty]}$ is said to be \emph{submultiplicative}, if there exists a constant $C>0$ such that  
$$
\Phi(ab) \le C \Phi(a) \Phi(b), \quad a,b>0.
$$
Similarly, $\Phi$ is said to be \emph{supermultiplicative} if there exists a constant $C>0$ such that
$$
\Phi(ab) \ge C \Phi(a) \Phi(b), \quad a,b>0.
$$
\end{defn}

\begin{example}
\begin{enumerate}
\item If $\Phi(t) = t + t^2$, then $\Phi$ is a submultiplicative Young function, since
$$
\Phi(ab) = ab + (ab)^2 \le ab+ab^2 + a^2b + (ab)^2 = \Phi(a)\Phi(b)
$$
for every $a,b\ge 0$.

\vrum

\item If $\Phi(t) = \begin{cases}
    t^3 & t<1, \\
    t & t\ge 1,
\end{cases}$
then $\Phi$ is a supermultiplicative Young function.

\vrum

\item For any $p\ge 1$ and $C>0$, $\Phi(t) = Ct^p$ is a Young function which is both sub- and supermultiplicative. In fact, in this case $\Phi(ab) = \frac{1}{C}\Phi(a)\Phi(b)$ for every $ab\ge 0$. By \autocite[Sec. 2.3, Prop. 12]{rao:1991}, up to equivalence, $\Phi(t) = t^p$, $1<p<\infty$ are the only N functions which are both sub- and supermultiplicative.
\end{enumerate}
\end{example}
% ----------------------------------------
%  Lebesgue exponents
% ----------------------------------------
\section{Orlicz spaces and Lebesgue exponents}\label{subsec:lebexp}
% ----------------------------------------------- %
\subsection{Lebesgue exponents}\label{subsec:expo}
% ----------------------------------------------- %

We begin by stating some observations about the properties of Young functions and some connections to their Lebesgue exponents $p_\Phi$ and $q_\Phi$. The following is a consequence of \autocite[Lemma 3.1]{mali:1985}.

\begin{prop}\label{prop:qpPhiIneq}
Let $\Phi\in\Delta_2$, $c\in[0,1]$, and $C\in[1,\infty)$. Then

\begin{equation}\label{eq:IneqqPhi}
c^{p_\Phi}\Phi(t) \le \Phi(ct) \le c^{q_\Phi} \Phi(t)
\end{equation}
and
\begin{equation}\label{eq:IneqpPhi}
C^{q_\Phi} \Phi(t) \le \Phi(Ct) \le C^{p_\Phi} \Phi(t).
\end{equation}
\end{prop}

\begin{proof}
For any $t>0$, by definition,
$$
q_\Phi \le \frac{t \Phi'(t)}{\Phi(t)} \le p_\Phi,
$$
or equivalently (since $\Phi(t) > 0$),
$$
t \Phi'(t)-p_\Phi \Phi(t) \le 0 \le  t \Phi'(t) - q_\Phi \Phi(t). 
$$
This implies that
$$
\frac{t^{p_\Phi}\Phi(t)- p_\Phi t^{p_\Phi - 1}\Phi'(t)}{t^{2p_\Phi}}\le 0 \le  \frac{t^{q_\Phi}\Phi(t) - q_\Phi t^{q_\Phi - 1}\Phi'(t)}{t^{2q_\Phi}}, 
$$
meaning
$$
\left( \frac{\Phi(t)}{t^{p_\Phi}}\right)' \le 0 \le \left( \frac{\Phi(t)}{t^{q_\Phi}}\right)'.
$$
Hence, there exists an increasing function $f(t)$ and a decreasing function $g(t)$ such that
$$
\Phi(t) = t^{q_\Phi} f(t) = t^{p_\Phi} g (t).
$$
If $c\in [0,1]$, then
$$
\Phi(c t) = (c t)^{q_\Phi}f(c t) \le c^{q_\Phi} t^{q_\Phi} f(t) = c^{q_\Phi} \Phi(t)
$$
and, using the other equality,
$$
\Phi(c t) = (c t)^{p_\Phi} g(c t) \ge c^{p_\Phi} \Phi(t)
$$
giving \eqref{eq:IneqqPhi}. The two inequalities of \eqref{eq:IneqpPhi} are obtained by arguing analogously.
\end{proof}

\begin{rem}
The assumption that $\Phi\in \Delta_2$ in \cref{prop:qpPhiIneq} is not necessary for the second inequality of \eqref{eq:IneqqPhi} or the first inequality of \eqref{eq:IneqpPhi} to hold; since $\Phi$ is a Young function, if $\Phi(s_1) = 0$ and $\Phi(s_2) = \infty$, then $s_1 < s_2$, $\Phi(t) = 0$ for all $t \le s_1$, and $\Phi(t) = \infty$ for all $t \ge s_2$. Repeating the arguments of the proof for $t\in (s_1,s_2)$ yields the results, since the aforementioned inequalities trivially hold for $t\in (0,s_1]\cup [s_2,\infty)$.
\end{rem}

\begin{rem}\label{rem:phidef}
Evidently, by the definition of $\nm{\cdot}{\Phi}$, we must have
$$
\rho_\Phi\left( \frac{f}{\nm{f}{\Phi}} \right) \le 1, \quad f\in L^\Phi,
$$
but if $p_\Phi < \infty$, then if $f\in L^\Phi$ and $\nm{f}{L^\Phi} \ne 0$, we must have
$$
\rho_\Phi\left(\frac{f}{\nm{f}{\Phi}} \right) = 1.
$$
To see why, suppose that
$$
\rho_\Phi \left( \frac{f}{\nm{f}{\Phi}}\right) = \alpha < 1.
$$
Let $\alpha^{1/p_\Phi} < c < 1$, and let $\lambda = c \nm{f}{\Phi}$. Then, by \cref{prop:qpPhiIneq},
$$
\rho_\Phi \left( \frac{f}{\lambda} \right) = \rho_\Phi \left( \frac{f}{c\nm{f}{\Phi}} \right) \le \frac{\alpha}{c^{p_\Phi}} < 1.
$$
But then $\lambda < \nm{f}{\Phi}$ and $\rho_\Phi(f/\lambda) < 1$, which contradicts the definition of $\nm{\cdot}{\Phi}$.
\end{rem}

\begin{rem}\label{rem:delta2}
A Young function is said to fulfill the \emph{$\Delta_2$-condition} \autocite[p. 6]{birn:1931} if there exists a constant $C>0$ such that
$$
\Phi(2t) \le C \Phi(t), \quad t>0.
$$
The fact that $p_\Phi < \infty$ implies $\Phi$ fulfills the $\Delta_2$-condition
is an easy consequence of the second inequality in \eqref{eq:IneqpPhi}. The converse is also easily shown since
$$
\Phi'(t) \le \frac{\Phi(2t) - \Phi(t)}{t} \le C\frac{\Phi(t)}{t},
$$
whereby $\Phi$ fulfilling the $\Delta_2$-condition for some constant $C>0$ implies that $\frac{t \Phi'(t)}{\Phi(t)}\le C$ for every $t>0$, hence $p_\Phi < \infty$. In the remainder of this paper, we write $\Phi\in\Delta_2$ to indicate that $\Phi$ is a Young function fulfilling the $\Delta_2$-condition.
\end{rem}

\begin{cor}\label{cor:pPhiqPhiineq}
Suppose that $\Phi\in\Delta_2$. Then
$$
\Phi(1) t^{p_\Phi} \le \Phi(t) \le \Phi(1) t^{q_\Phi}, \quad t\le 1
$$
and
$$
\Phi(1) t^{q_\Phi} \le \Phi(t) \le \Phi(1) t^{p_\Phi}, \quad t\ge 1.
$$
\end{cor}

\begin{proof}
Apply \eqref{eq:IneqqPhi} and \eqref{eq:IneqpPhi} with $t=1$.
\end{proof}

We can extend this to the following.

\begin{cor}\label{cor:normPsiprop}
Suppose that $\Phi\in\Delta_2$ and let $\Psi(t) = \frac{\Phi(t)}{\Phi(1)}$. Then the following holds true.
\begin{enumerate}
\item \label{cor:propitem1}$L^\Phi(\rd) = L^\Psi(\rd)$;

\vrum

\item \label{cor:propitem2}$p_\Phi = p_\Psi$ and $q_\Phi=q_\Psi$;

\vrum

\item \label{cor:propitem3}
$$
t^{p_\Psi} \le \Psi(t) \le t^{q_\Psi}, \qquad t\le 1
$$ 
and
$$
t^{q_\Psi} \le \Psi(t) \le t^{p_\Psi}, \qquad t > 1;
$$

\vrum

\item \label{cor:propitem4}
$$
q_\Psi t^{p_\Psi-1} \le \Psi'(t) \le p_\Psi t^{q_\Psi-1}, \qquad t \le 1,
$$
and
$$
q_\Psi t^{q_\Psi-1} \le \Psi'(t) \le p_\Psi t^{p_\Psi-1}, \qquad t > 1;
$$

\vrum

\item \label{cor:propitem5}the inverse $\Psi^{-1}$ exists and fulfills 
$$
t^{1/q_\Psi} \le \Psi^{-1}(t) \le t^{1/p_\Psi}, 
\qquad t\le 1,
$$
and
$$
t^{1/p_\Psi} \le \Psi^{-1}(t) \le t^{1/q_\Psi}, 
\qquad t\ge 1.
$$
\end{enumerate}
\end{cor}

\begin{proof}
Statements \ref{cor:propitem1}-\ref{cor:propitem3} follow directly from definitions and \cref{cor:pPhiqPhiineq}.

For statement \ref{cor:propitem4}, we have 
$$
q_\Psi \Psi(t)\le t \Psi'(t) \le p_\Phi \Psi(t)
$$
and hence the statement follows by \ref{cor:propitem3}.
Finally, \ref{cor:propitem5} follows from \ref{cor:propitem3} by noting that
$t^{q_\Phi} \le \Phi(t)$ implies $\Phi^{-1}(t^{q_\Phi}) \le t$, so that, with $s=t^{q_\Phi}$, $\Phi^{-1}(s) \le s^{1/q_\Phi}$. Since $s\le 1$ whenever $t\le 1$, the first inequality in \ref{cor:propitem5} holds. The others follow by analogous arguments.
\end{proof}

We now present some examples of Young functions $\Phi$ with different Lebesgue exponents $q_\Phi$ and $p_\Phi$. The values of these Lebesgue exponents in \ref{ex:item1} and \ref{ex:item2} follow by very simple calculations. The calculations for \ref{ex:item5} and \ref{ex:item6} are left for the reader.

\begin{example}\label{ex:pPhiqPhi} We have the following.
\begin{enumerate}
\item \label{ex:item1} Let $\Phi(t) = t^p$, $p\ge 1$. Then $p_\Phi = q_\Phi = p$.

\vrum

\item \label{ex:item2} Let $\Phi(t) = t^p + t^q$, $ p > q \ge 1$. Then
$p_\Phi = p$, $q_\Phi = q$.

\vrum

\item\label{ex:item3} Let
$$
\Phi(t) = 
\begin{cases}
\frac12 t^{2}, & t\le 1 \\
t - \frac12 & 1< t \le 2 \\
\frac14 t^{2} + \frac12 & t>2.
\end{cases}
$$
Then, since $\Phi'_+(1) = \Phi'_-(1)$, $\Phi'_+(2) = \Phi'_-(2)$ and $\Phi''(t)\ge 0$ for all $t>0$, $\Phi$ is a Young function and
$$
\frac{t \Phi'(t)}{\Phi(t)} =
\begin{cases}
2, & t\le 1 \\
\dfrac{1}{1 - \frac{1}{2 t}}, & 1< t \le 2 \\[10pt]
\dfrac{2}{1 + \frac{2}{t^2}} & t>2.
\end{cases}
$$
Hence 
\begin{align*}
q_\Phi 
&= \min\left\{2, 
\inf_{1< t \le 2}\left(\frac{1}{1 - \frac{1}{2 t}}\right), 
\inf_{t>2}\left(\frac{2}{1 + \frac{2}{t^2}}\right)\right\} \\
&=\frac43,
\end{align*}
and
\begin{align*}
p_\Phi 
&= \max\left\{2, 
\sup_{1< t \le 2}\left(\frac{1}{1 - \frac{1}{2 t}}\right), 
\sup_{t>2}\left(\frac{2}{1 + \frac{2}{t^2}}\right)\right\} \\
&=2
\end{align*}
\item\label{ex:item4} Let
$$
\Phi(t) = 
\begin{cases}
\frac32 t^{2}, & t\le 1 \\
t^3 + \frac12 & 1< t \le 2 \\
3 t^{2} - \frac72 & t>2.
\end{cases}
$$
Then, $\Phi$ is a Young function, 
\begin{align*}
q_\Phi 
&= \min\left\{2, 
\inf_{1< t \le 2}\left(\frac{3}{1 + \frac{1}{2 t^3}}\right), 
\inf_{t>2}\left(\frac{6}{3 - \frac{7}{2t^2}}\right)\right\} \\
&=2,
\end{align*}
and
\begin{align*}
p_\Phi 
&= \max\left\{2, 
\sup_{1< t \le 2}\left(\frac{3}{1 + \frac{1}{2 t^3}}\right), 
\sup_{t>2}\left(\frac{6}{3 - \frac{7}{2t^2}}\right)\right\} \\
&=\frac{48}{17}.
\end{align*}

\vrum

\item \label{ex:item5} Let $\Phi(t) = t^n \ln^m (1 + t)$, where $n\in \na\setminus\{0\}$, $m\in \na$. Then $q_\Phi = n$, $p_\Phi = n + m$.

\vrum

\item \label{ex:item6} Let $\Phi(t) = t^n e^t$, $n\in\na \setminus\{0\}$. Then $q_\Phi = n$, $p_\Phi = \infty$.
\end{enumerate}
\end{example}

\begin{rem}
The Young functions of \cref{ex:pPhiqPhi}\ref{ex:item3}-\ref{ex:item4} highlight the limitations of analysing Orlicz spaces by way of Lebesgue exponents, since both of these Young functions are equivalent to $\Psi(t) = t^{2}$, despite having different Lebesgue exponents. We can expand this type of construction further.
\end{rem}

\begin{example}
Let $r_1,r_2 \ge 1$, $r_1\neq r_2$ and let
$$
\Phi(t) = 
\begin{cases}
\frac{c_1}{r_1} t^{r_1}, & t\le a, \\[5pt]
\frac{1}{r_2}t^{r_2} + d_1, & a< t \le b, \\[5pt]
\frac{c_2}{r_1} t^{r_1} + d_2, & t>b,
\end{cases}
$$
where 
\begin{align*}
c_1 &= a^{r_2-r_1}, &
d_1 &= \left( \frac{1}{r_1} - \frac{1}{r_2}\right)a^{r_2}, \\
c_2 &= b^{r_2 - r_1}, \text{ and} &
d_2 &= \left(\frac{1}{r_1} - \frac{1}{r_2}\right)\left( a^{r_2} - b^{r_2}\right).
\end{align*}
Then $\Phi$ is a Young function which is equivalent to $\Psi(t) = t^{r_1}$. On the other hand,
$$
\frac{t \Phi'(t)}{\Phi(t)} =
\begin{cases}
r_1, & t\le a, \\
\dfrac{r_2}{1 + k_1 t^{-r_2}}, & a< t \le b, \\[10pt]
\dfrac{r_1}{1 + k_2 t^{-r_1}}, & t>b,
\end{cases}
$$
where 
$$
k_1 = a^{r_2}\left(\frac{r_2}{r_1} - 1\right), \text{ and}\quad
k_2 = b^{r_1}\left( 1 - \left(\frac{a}{b}\right)^{r_2}  \right) \left( \frac{r_1}{r_2} - 1\right).
$$
If $r_2 > r_1$, then $k_1> 0$ and $k_2<0$ so 
\begin{align*}
q_\Phi 
&= \min
\left\{r_1, 
\inf_{a< t \le b}\frac{r_2}{1+k_1 t^{-r_2}}, 
\inf_{t>b}\frac{r_1}{1 + k_2 t^{-r_1}}
\right\} \\
&=\min\{r_1, \frac{r_2}{1+k_1 a^{-r_2}},r_1\} = r_1,
\end{align*}
and
\begin{align*}
p_\Phi 
&= \max
\left\{r_1, 
\sup_{a< t \le b}\frac{r_2}{1+k_1 t^{-r_2}}, 
\sup_{t>b}\frac{r_1}{1 + k_2 t^{-r_1}}
\right\} \\
&=\max\left\{r_1,\frac{r_2}{1 + k_1 b^{-r_2}},\frac{r_1}{1+k_2 b^{-r_1}}\right\} \le r_2,
\end{align*}
and, if we let $b\rightarrow \infty$ then $p_\Phi \rightarrow r_2$.
On the other hand, if $r_1 > r_2$, then $k_1<0$ and $k_2>0$ so
$$
q_\Phi = \min\left\{r_1,\frac{r_2}{1 + k_1 b^{-r_2}},\frac{r_1}{1+k_2 b^{-r_1}}\right\}\ge r_2,
$$
$$
p_\Phi =\sup\left\{r_1, \frac{r_2}{1+k_1 a^{-r_2}},r_1\right\} = r_1,
$$  
and in this case, $b\rightarrow \infty$ implies $q_\Phi \rightarrow r_2$.
\end{example}

From this example, we can conclude the following.

\begin{prop}
Let $1< p_1 < p < p_2 <\infty$. Then there exist Young functions $\Phi$, $\Psi$ with $q_\Phi < p_1$, $p_\Psi > p_2$, and $L^\Phi = L^\Psi = L^p$.
\end{prop}

\begin{proof}
Let $\Phi$ be as in the previous example with $1<r_2<p_1$ and $r_1=p$. Since $q_\Phi \rightarrow r_2$ as $b\rightarrow \infty$, we simply need to choose $b$ large enough and the result holds. By choosing $\Psi$ in the same way but with $r_2>p_2$ instead, we obtain $p_\Psi > p_2$ in the same way.
\end{proof}

Of course, we can simply add another finite segment where the function has the form $\frac{c_3}{r_3}t^{r_3} + d_3$. Refining the arguments of the previous example, we obtain the following proposition.

\begin{prop}\label{prop:p1p2Lp}
Let $1 < p_1 < p < p_2 < \infty$. Then there exists a Young function $\Phi$ with $L^\Phi = L^p$, $q_\Phi = p_1$, and $p_\Phi = p_2$.
\end{prop}

\begin{proof}
Let $1<\alpha < \beta$, $1<r_1 < p_1$, $r_2 > p_2$, and
$$
\Phi_{\alpha,\beta} (t) =
\begin{cases}
\frac{t^p}{p}, & t\le 1, \\[5pt]
\frac{c_1 t^{r_1}}{r_1} + d_1, & 1 < t \le \alpha, \\[5pt]
\frac{c_2 t^{r_2}}{r_2} + d_2, & \alpha < t \le \beta, \\[5pt]
\frac{c_3 t^p}{p} + d_3, & t>\beta.
\end{cases}
$$
To ensure that 
$$
\lim_{t\rightarrow t_j^+} \Phi_{\alpha,\beta}'(t) 
= \lim_{t\rightarrow t_j^-} \Phi_{\alpha,\beta}'(t), \qquad
\lim_{t\rightarrow t_j^+} \Phi_{\alpha,\beta}(t) 
= \lim_{t\rightarrow t_j^-} \Phi_{\alpha,\beta}(t)
$$ for $t_1 = 1, t_2 = \alpha,$ and $ t_3 = \beta$, we set
\begin{align*}
c_1 &= 1, &
c_2 &= \alpha^{r_1 - r_2}, & 
c_3 &= \alpha^{r_1 - r_2} \beta^{r_2 - p}, \\
d_1 &= \frac{1}{p} - \frac{1}{r_1}, &
d_2 &= \alpha^{r_1}\left(\frac{1}{r_1} - \frac{1}{r_2} \right) + d_1,&
d_3 &= \alpha^{r_1-r_2}\beta^{r_2}\left(\frac{1}{r_2} - \frac{1}{p} \right) + d_2.
\end{align*}
It is clear that, for fixed $\alpha,\beta$ fulfilling $1<\alpha<\beta$, $\Phi_{\alpha,\beta}(t)$ is equivalent to $t^p$ for $t\le 1$ and $t>\beta$, and since $(\alpha,\beta]$ is a finite interval, we simply need to choose $C>0$ sufficiently large to obtain $C^{-1} \Phi_{\alpha,\beta}(t) \le t^p \le C\Phi_{\alpha,\beta}(t)$ for all $t\ge 0$. Hence $L^{\Phi_{\alpha,\beta}} = L^p$ for any choice of $\beta>\alpha>1$.
Since
$$
\Phi_{\alpha,\beta}'(t) =
\begin{cases}
t^{p- 1}, & t\le 1, \\
t^{r_1-1}, & 1 < t \le \alpha, \\
c_2 t^{r_2- 1}, & \alpha < t \le \beta, \\
c_3 t^{p-1}, & t>\beta,
\end{cases}
$$
we have
$$
\frac{t \Phi_{\alpha,\beta}'(t)}{\Phi_{\alpha,\beta} (t)}=
\begin{cases}
p, & t\le 1, \\
\dfrac{r_1}{1 + k t^{-r_1}}, & 1 < t \le \alpha, \\[10pt]
\dfrac{r_2}{1+ l_\alpha t^{-r_2}}, & \alpha < t \le \beta, \\[10pt]
\dfrac{p}{1+ m_{\alpha,\beta} t^{-p}}, & t>\beta,
\end{cases}
$$
where
\begin{align*}
k &=r_1 d_1 =r_1\left(\frac{1}{p} - \frac{1}{r_1}\right), \\ 
l_\alpha &= \frac{r_2 d_2}{c_2} 
= r_2 \alpha^{r_2}\left[\left(\frac{1}{r_1} - \frac{1}{r_2} \right) 
+ \alpha^{-r_1}d_1\right] = r_2 \alpha^{r_2} C(\alpha),\\
m_{\alpha,\beta} &= \frac{p d_3}{c_3} =  p \beta^p \left[ \left(\frac{1}{r_2} - \frac{1}{p} \right)
+ \left(\frac{\alpha}{\beta}\right)^{r_2}C(\alpha)\right].
\end{align*}
Hence
\begin{align*}
&q_{\Phi_{\alpha,\beta}} \\
&= 
\min\left\{
p, 
\inf_{1<t\le \alpha} \left( \frac{r_1}{1 + k t^{-r_1}}\right),
\inf_{\alpha<t\le \beta} \left( \frac{r_2}{1+ l_\alpha t^{-r_2}}\right),
\inf_{t\ge \beta} \left( \frac{p}{1+ m_{\alpha,\beta} t^{-p}}\right) 
\right\}. 
\end{align*}
It is clear that $k<0$, $l_\alpha > 0$ for all $\alpha>1$, but that the sign of $m_{\alpha,\beta}$ can vary depending on the choices of $\alpha$ and $\beta$. If $m_{\alpha,\beta} < 0$ then
$$
p < \frac{p}{1+ m_{\alpha,\beta} t^{-p}} \le \frac{p}{1+ m_{\alpha,\beta} \beta^{-p}} = \frac{r_2}{1+ l_\alpha \beta^{-r_2}}
$$
and if $m_{\alpha,\beta} > 0$ then
$$
\frac{r_2}{1+ l_\alpha \beta^{-r_2}} \le \frac{p}{1+ m_{\alpha,\beta} t^{-p}} < p
$$
for $t \ge \beta$. In either case, neither the infimum nor the supremum of $\frac{t \Phi_{\alpha,\beta}'(t)}{\Phi_{\alpha,\beta} (t)}$ will be determined by its value for $t>\beta$. Therefore
\begin{align*}
q_{\Phi_{\alpha,\beta}} &= \min
\left\{ 
p, 
\frac{r_1}{1 + k \alpha^{-r_1}}, 
\frac{r_2}{1+ l_\alpha \alpha^{-r_2}}
\right\} \\
&=\frac{r_1}{1 + k \alpha^{-r_1}}.
\end{align*}
Let $h(\alpha)=\frac{r_1}{1 + k \alpha^{-r_1}}$. Then $h(\alpha)$ depends continuously on $\alpha$,
$$
h(1) = \frac{r_1}{1+k}, \quad\text{and}\quad h(\alpha)\rightarrow r_1, \quad \alpha\rightarrow \infty.
$$
Since
$$
r_1 < p_1 < \frac{r_1}{1+k} = p,
$$ 
there must exist $\gamma>1$ such that 
$$
q_{\Phi_{\gamma,\beta}} = h(\gamma) = \frac{r_1}{1 + k \gamma^{-r_1}} = p_1
$$
for all $\beta > \gamma$.  Furthermore,
\begin{align*}
&p_{\Phi_{\gamma,\beta}} \\
&= 
\max\left\{
p, 
\sup_{1<t\le \gamma} \left( \frac{r_1}{1 + k t^{-r_1}}\right),
\sup_{\gamma<t\le \beta} \left( \frac{r_2}{1+ l_{\gamma} t^{-r_2}}\right),
\sup_{t\ge \beta} \left( \frac{p}{1+ m_{\gamma,\beta} t^{-p}}\right) 
\right\} \\
&= 
\max
\left\{ 
p, 
\frac{r_1}{1 + k}, 
\frac{r_2}{1+ l_{\gamma}  \beta^{-r_2}}
\right\} \\
&= \frac{r_2}{1+ l_{\gamma}\beta^{-r_2}} \rightarrow r_2, \quad
\beta\rightarrow \infty.
\end{align*}
Since
$$
r_2 > \frac{r_2}{1+ l_{\gamma}\beta^{-r_2}} \ge  \frac{r_2}{1+ l_{\gamma}\gamma^{-r_2}} = \frac{r_1}{1 + k \gamma^{-r_1}} = p_1
$$
and $p_1 < p_2 < r_2$, there similarly exists a $\delta> \gamma$ such that 
$$
p_{\Phi_{\gamma,\delta}} = \frac{r_2}{1+ l_{\gamma}\delta^{-r_2}} = p_2.
$$ 
With $\Phi = \Phi_{\gamma,\delta}$, the desired result has now been reached.
\end{proof}

By slightly modifying this type of construction, we obtain the following.

\begin{prop}\label{prop:p1p2LPhi}
Let $\Phi$ be a Young function and let $1 < p_1 < q_\Phi \le p_\Phi < p_2 < \infty$. Then there exists a Young function $\Psi$ with $L^\Psi = L^\Phi$, $q_\Psi \le p_1$, and $p_\Psi \ge p_2$.
\end{prop}

\begin{proof}
The proof is similar to that of \cref{prop:p1p2Lp}.
Let $1<\alpha < \beta$, $1<r_1 < p_1$, $r_2 > p_2$, and
$$
\Psi_{\alpha,\beta} (t) =
\begin{cases}
\frac{\Phi(t)}{\Phi'(1)} & t\le 1, \\[5pt]
\frac{c_1 t^{r_1}}{r_1} + d_1 & 1 < t \le \alpha \\[5pt]
\frac{c_2 t^{r_2}}{r_2} + d_2 & \alpha < t \le \beta \\[5pt]
\frac{c_3 \Phi(t)}{\Phi'(\beta)} + d_3 & t>\beta.
\end{cases}
$$
In this case, with $g_\Phi(t) = \frac{t\Phi'(t)}{\Phi(t)}$, we set
\begin{align*}
c_1 &= 1, &
c_2 &= \alpha^{r_1 - r_2}, & 
c_3 &= \alpha^{r_1 - r_2} \beta^{r_2 - 1}, \\
d_1 &= \frac{1}{g_\Phi(1)} - \frac{1}{r_1}, &
d_2 &= \alpha^{r_1}\left(\frac{1}{r_1} - \frac{1}{r_2} \right) + d_1,&
d_3 &= \alpha^{r_1-r_2}\beta^{r_2}\left(\frac{1}{r_2} - \frac{1}{g_\Phi(\beta)} \right) + d_2.
\end{align*}
Continuing as in the proof of \cref{prop:p1p2Lp}, we obtain
$$
\frac{t \Psi_{\alpha,\beta}'(t)}{\Psi_{\alpha,\beta} (t)}=
\begin{cases}
g_\Phi(t) & t\le 1 \\
\dfrac{r_1}{1 + k t^{-r_1}} & 1 < t \le \alpha \\[10pt]
\dfrac{r_2}{1+ l_\alpha t^{-r_2}} & \alpha < t \le \beta \\[10pt]
\dfrac{g_\Phi(t)}{1 + m_{\alpha,\beta}(\Phi(t))^{-1}} & t>\beta.
\end{cases}
$$
where
\begin{align*}
k &=r_1 d_1 
= r_1 \left( \frac{1}{g_\Phi(1)} - \frac{1}{r_1}\right), 
\\ 
l_\alpha &= \frac{r_2 d_2}{c_2} 
= r_2\alpha^{r_2}\left[\left(\frac{1}{r_1} - \frac{1}{r_2} \right) 
+ \alpha^{-r_1} d_1\right] = r_2 \alpha^{r_2}C(\alpha),
\\
m_{\alpha,\beta} &= \frac{\Phi'(\beta) d_3}{c_3} 
= \beta \Phi'(\beta)\left[\left( \frac{1}{r_2} - \frac{1}{g_\Phi(\beta)}\right) 
+ \left(\frac{\alpha}{\beta}\right)^{r_2} C(\alpha)\right].
\end{align*}
As before, we have $k<0$ and $l_\alpha > 0 $ for all $\alpha >1$.

In this case, depending on which values are attained by $g_\Phi(t)$, it is possible that the infimum or supremum is attained for $t>\beta$, and that this is less than $r_1$ or greater than $r_2$. For sufficiently large $\alpha$, we have $m_{\alpha,\beta}< 0$ and the arguments proceed as in the previous proof. But then the infimum may (depending on the value of $q_\Phi/p_\Phi$) be smaller than $p_1$ for all such $\alpha,\beta$. For this reason, it is difficult to find a specific $\gamma>1$ such that $q_{\Psi_{\gamma,\beta}} = p_1$. However, it is the case that
$$
q_{\Psi_{\alpha,\beta}} \le \inf_{1< t\le \alpha} \left(\dfrac{r_1}{1 + k t^{-r_1}}\right) = \frac{r_1}{1 + k\alpha^{-r_1}} \rightarrow r_1 < p_1,
$$
hence there exists $\gamma > 1$ such that $q_{\Psi_{\gamma,\beta}} < p_1$ for all $\beta > \gamma$.

Similarly, we obtain
$$
p_{\Psi_{\gamma,\beta}} \ge \sup_{\alpha < t\le \beta}\left(\dfrac{r_2}{1+ l_\gamma t^{-r_2}}\right) = \dfrac{r_2}{1+ l_\gamma \beta^{-r_2}} \rightarrow r_2 > p_2,
$$
hence there exists $\delta > \gamma$ such that $p_{\Psi_{\gamma,\delta}} > p_2$. This completes the proof.
\end{proof}

%%

% ----------------------------------------------------------- %
\subsection{Equivalent Young functions}\label{subsec:equivYfn}
% ----------------------------------------------------------- %

The preceding proposition poses a difficult problem when attempting to analyze Orlicz spaces based on the Lebesgue exponents of the corresponding Young functions. In an attempt to avoid this problem, we will consider equivalence classes of Young functions and define the Lebesgue exponents $q_\Phi$ and $p_\Phi$ on these classes as a whole. This provides an opportunity to improve certain results, although there are limitations to this approach as well.

We note that if $\Phi$ and $\Psi$ are equivalent Young functions, then their (semi) derivatives are also equivalent, although they may not be Young functions.

\begin{prop}
Suppose $\Phi,\Psi$ are equivalent Young functions and $p_\Phi,p_\Psi < \infty$. Then there exists a constant $C_1 \ge 1$ such that
$$
C_1^{-1} \Psi' \le \Phi' \le C_1 \Psi'. 
$$
\end{prop}

\begin{proof}
Let $C\ge 1$ be a constant fulfilling
$$
C^{-1} \Psi \le \Phi \le C\Psi.
$$
Since
$$
q_\Phi \le \frac{t \Phi'(t)}{\Phi(t)} \le p_\Phi, \quad q_\Psi \le \frac{t\Psi'(t)}{\Psi(t)} \le p_\Psi, \quad t>0,
$$
we obtain
$$
\Phi'(t) \le p_\Phi \frac{\Phi(t)}{t} \le C p_\Phi \frac{\Psi(t)}{t} \le \frac{C p_\Phi}{q_\Psi} \Psi'(t),
$$
and similarly
$$
\Psi'(t) \le p_\Psi \frac{\Psi(t)}{t} \le C p_\Psi \frac{\Phi(t)}{t} \le \frac{C p_\Psi}{q_\Phi} \Phi'(t),
$$
hence the desired inequality is fulfilled with
\begin{equation*}
C_1 = \max\left\{ \frac{C p_\Phi}{q_\Psi}, \frac{C p_\Psi}{q_\Phi}\right\}. \qedhere
\end{equation*}
\end{proof}

As \cref{prop:p1p2LPhi} shows, $\Phi$ and $\Psi$ can be equivalent even if $p_\Phi \neq p_\Psi$, $q_\Phi \neq q_\Psi$. Conversely, $p_\Phi = p_\Psi$, $q_\Phi = q_\Psi$ does not imply that the Young functions are equivalent. This can be seen by considering, for example, 
$$
\Phi(t) = t^2 \ln(1+t) \quad \text{and }\quad \Psi(t) = t^2 + t^3.
$$ 
Then $p_\Phi = p_\Psi = 3$ and $q_\Phi = q_\Psi = 2$, but $\Phi$ and $\Psi$ are not equivalent (since, for instance, with $d=1$ and $f(x) = (1+x^2)^{-1/4}$, $\nm{f}{\Phi} < \infty$, but $\nm{f}{\Psi} = \infty$.)

Since equivalent Young functions have equivalent growth rates, we define Lebesgue exponents for equivalent classes as follows.

\begin{defn}
Let $\Phi$ be a Young function and $[\Phi]$ its equivalence class. We define
$$
p_{[\Phi]} = \inf\Sets{p_\Psi}{\Psi\in [\Phi]}
$$
and
$$
q_{[\Phi]} = \sup\Sets{q_\Psi}{\Psi\in [\Phi]}.
$$
\end{defn}

The following proposition follows immediately from \autocite[Corollary 2.8]{boni:2024}.

\begin{prop}
Suppose that $\Phi$ is a Young function. Then the following statements are equivalent:
\begin{enumerate}

\item $p_{[\Phi]} < \infty$ $(q_{[\Phi]} > 1)$;

\vrum

\item $p_\Psi < \infty$ $(q_\Psi > 1)$ for every $\Psi\in[\Phi]$;

\vrum

\item $p_\Psi < \infty$ $(q_\Psi > 1)$ for some $\Psi\in[\Phi]$.

\end{enumerate}
\end{prop}

This means that if $p_\Phi = \infty$ and $q_\Phi = 1$, then the same is true for all $\Psi$ in the equivalence class $[\Phi]$, and further $p_{[\Phi]} = \infty$ and $q_{[\Phi]} = 1$. However, this does not mean that any $\Phi_1$, $\Phi_2$ with $p_{\Phi_j} = \infty$, $q_{\Phi_j} = 1$, $j=1,2$ are equivalent. One counterexample is given by
$$
\Phi_1(t) = \begin{cases}
0, & t=0, \\
e^{-1/t}, & 0<t<1/2, \\
e^{-2}(4t - 1), & t\ge 1/2, 
\end{cases}
\mathandd
\Phi_2(t) = e^{t}  - 1.
$$

We now provide some examples of operations which preserve equivalence classes of Young functions. The details are left for the reader. For more information on the second of the following three examples (in the case that the Young functions involved are not equivalent), as well as several other examples, see \autocite[p. 98-99]{mali:1989B}.

\begin{example}
Suppose that $\Phi_1,\Phi_2,\Psi$ are Young functions and that $\Phi_1,\Phi_2 \in [\Psi]$.
\begin{enumerate}
\item If $\Phi = \alpha \Phi_1 + \beta\Phi_2$, $\alpha,\beta>0$, then $\Phi \in [\Psi]$, 
$$
p_\Phi 
\le 
\min\left\{\frac{C}{C+1}(p_{\Phi_1} + p_{\Phi_2}),\max\{p_{\Phi_1},p_{\Phi_2}\}\right\}
$$ and 
$$
q_\Phi \ge \max\left\{\frac{1}{1 + C}(q_{\Phi_1} + q_{\Phi_2}), \min\{q_{\Phi_1},q_{\Phi_2}\}\right\}
$$
for some constant $C\ge 1$.

\vrum

\item If $\Phi = \max\{\Phi_1,\Phi_2\}$, then $\Phi\in[\Psi]$,
$$
p_\Phi \le \max\{p_{\Phi_1},p_{\Phi_2}\}, \quad 
q_\Phi\ge  \min\{q_{\Phi_1},q_{\Phi_2}\}
$$

\vrum

\item If
$$
\Phi = \sup \Sets{\Phi_3 \text{ convex}}{\Phi_3 \le \min\{\Phi_1,\Phi_2\}},
$$
then $\Phi$ is a Young function and $\Phi \in [\Psi]$.
\end{enumerate}
\end{example}

It will become clear that $p_{[\Phi]}$ and $q_{[\Phi]}$ above are closely related to how the function
$$
g_\Phi(t) = \frac{t \Phi'(t)}{\Phi(t)}
$$
behaves at infinity and at the origin. The quantities $\limsup_{t\rightarrow \omega} g_\Phi(t)$ and $\liminf_{t\rightarrow \omega} g_\Phi(t)$, where $\omega=0$ or $\omega=\infty$, appear in the literature (cf. e.g. \autocite[Theorem 11.11]{mali:1989B} or \autocite[Corollary 4]{rao:1991}). In the case that the limit of $g_\Phi(t)$ exists as $t$ approaches infinity or the origin, one can say the following.

\begin{lem}\label{lem:equivp1p2}
 Suppose that $\Phi$ and $\Psi$ are equivalent Young functions. 
 \begin{enumerate}
\item \label{lem:equivitem1} If
$$
\lim_{t\rightarrow \infty} \frac{t \Phi'(t)}{\Phi(t)} = p_1 \quad \text{and}\quad \lim_{t\rightarrow \infty} \frac{t \Psi'(t)}{\Psi(t)} = p_2,
$$
then $p_1 = p_2$.
\item \label{lem:equivitem2} If
$$
\lim_{t\rightarrow 0} \frac{t \Phi'(t)}{\Phi(t)} = q_1 \quad \text{and}\quad \lim_{t\rightarrow 0} \frac{t \Psi'(t)}{\Psi(t)} = q_2,
$$
then $q_1 = q_2$.
\end{enumerate}
\end{lem}

\begin{proof}
Suppose that $p_1 < p_2$.
For every $\varepsilon > 0$, there exists a $K_\varepsilon>0$ such that
$$
- \varepsilon < \frac{t\Phi'(t)}{\Phi(t)}  - p_1 < \varepsilon, \qquad t>K_\varepsilon,
$$
which implies
$$
\left(t^{-(p_1+\varepsilon)} \Phi(t)\right)' < 0, \qquad \left( t^{-(p_1 - \varepsilon)}\Phi(t) \right)' > 0 
$$
hence
$$
\Phi(t) = t^{p_1 + \varepsilon} f(t) = t^{p_2 - \varepsilon} g(t),
$$
where $f(t) = t^{-(p_1+\varepsilon)}\Phi(t)$ is a decreasing function, and $g(t) =t^{-(p_1 - \varepsilon)}\Phi(t)$ is an increasing function.
Hence
$$
A_\varepsilon t^{p_1 - \varepsilon} \le \Phi(t) \le B_\varepsilon t^{p_1 + \varepsilon}, \quad t> K_\varepsilon,
$$
where $A_\varepsilon = g(K_\varepsilon)$ and $B_\varepsilon = f(K_\varepsilon)$. By a similar argument, there exist $C_\varepsilon, D_\varepsilon, L_\varepsilon > 0$ such that
$$
C_\varepsilon t^{p_2 - \varepsilon} \le \Psi(t) \le D_\varepsilon t^{p_2 + \varepsilon}, \quad t> L_\varepsilon.
$$
By picking $\varepsilon < \frac{p_2 - p_1}{2}$, we get for $t> \max\{K_\varepsilon,L_\varepsilon \}$
$$
t^{p_1 + \varepsilon} < t^{p_2 - \varepsilon} \le C_\varepsilon^{-1} \Psi(t) \le C C_\varepsilon^{-1}\Phi(t) \le \frac{CB_\varepsilon}{C_\varepsilon} t^{p_1 + \varepsilon}.
$$
This implies that, for $C_1 = \frac{C B_\varepsilon}{C_\varepsilon}$ and $t>\max\{K_\varepsilon, L_\varepsilon\}$,
$$
C_1^{-1}t^{p_1 + \varepsilon} \le t^{p_2 - \varepsilon} \le C_1 t^{p_1 + \varepsilon}
$$
which is a contradiction. Hence \ref{lem:equivitem1} holds, and \ref{lem:equivitem2} follows by analogous computations.
\end{proof}

As before, let $g_\Phi(t) = \frac{t \Phi'(t)}{\Phi(t)}$. If $g_\Phi(t)$ approaches $p_0$ and $p_\infty$ as $t$ approaches the origin and infinity, respectively, then it is clear that
$$
q_\Phi \le p_0,p_\infty \le p_\Phi.
$$
We can, however, relate these quantities more closely.

\begin{prop}\label{prop:p1p2pPhiqPhi}
Suppose that $\Phi$ is a Young function,
$$
p_0 =\lim_{t\rightarrow 0} \frac{t \Phi'(t)}{\Phi(t)}, \mathand 
p_\infty = \lim_{t\rightarrow \infty} \frac{t \Phi'(t)}{\Phi(t)}.
$$
Let $p = \max\{p_0,p_\infty\}$ and $q=\min\{p_0,p_\infty\}$.
Then, for every $\varepsilon > 0$, there exists a Young function $\Psi \in [\Phi]$ with 
$$
p\le p_\Psi < p+\varepsilon \mathand 
q-\varepsilon < q_\Psi \le q.
$$

\end{prop}

\begin{proof}
Suppose $\varepsilon > 0$. Let $g_\Phi(t) = \frac{t \Phi'(t)}{\Phi(t)}$.
Suppose $q \le r \le p$ and for $n>1$ let
$$
\Psi_n(t) = 
\begin{cases}
\frac{\Phi(t)}{\Phi'(1/n)}, &t\le 1/n, \\[5pt]
\frac{a t^r}{r} + b, & 1/n < t \le n, \\[5pt]
\frac{c \Phi(t)}{\Phi'(n)} + d, & t>n,
\end{cases}
$$
where
\begin{align*}
a &= n^{r-1}, &
b &= n^{-1}\left(\frac{1}{g_\Phi(1/n)} - \frac{1}{r}\right), \\
c &= n^{2r - 2}, &
d &= n^{2r-1}\left(\frac{1}{r} - \frac{1}{g_\Phi(n)}\right)+ n^{-1}\left(\frac{1}{g_\Phi(1/n)} - \frac{1}{r}\right). 
\end{align*}
We have 
$$
\Psi'_n(t) = 
\begin{cases}
\frac{\Phi'(t)}{\Phi'(1/n)}, &t\le 1/n, \\[5pt]
a t^{r-1}, & 1/n < t \le n, \\[5pt]
\frac{c \Phi'(t)}{\Phi'(n)}, & t>n,
\end{cases}
$$
hence
$$
\frac{t\Psi_n' (t)}{\Psi_n(t)} =
\begin{cases}
g_\Phi(t), & t\le 1/n, \\[5pt]
\dfrac{r}{1 + k_nt^{-r}}, & 1/n < t \le n, \\[10pt]
\dfrac{g_\Phi(t)}{1 + l_n \cdot (\Phi(t))^{-1}}, & t>n,
\end{cases}
$$
where
$$
k_n = \frac{r b}{a} = n^{-r}\left( \frac{r}{g_\Phi(1/n)} - 1 \right)
$$
and
$$
l_n = \frac{\Phi'(n) d}{c} = \Phi(n) \left[ \left( \frac{g_\Phi(n)}{r} - 1\right) + n^{-2 r} \left( \frac{g_\Phi(n)}{g_\Phi(1/n)} - \frac{g_\Phi(n)}{r}\right) \right].
$$
Since 
$$
\frac{t\Psi_n'(t)}{\Psi_n(t)} = g_\Phi(t), \qquad t<1/n,
$$
it is clear that $p_{\Psi_n} \ge p_1 \ge q_{\Psi_n}$, and since 
$$
\frac{g_\Phi(t)}{1 + l_n \cdot(\Phi(t))^{-1}} \rightarrow p_2, \qquad t\rightarrow \infty,
$$
we must have $p_{\Psi_n} \ge p_2 \ge q_{\Psi_n}$ as well. Hence $p_{\Psi_n} \ge p$ and $q_{\Psi_n} \le q$.

Let $0<\varepsilon_1 < \varepsilon$ be chosen such that $\frac{p_2 - \varepsilon_1}{p_2+\varepsilon_1} >  1 - \frac{\varepsilon}{q}$ and $\frac{p_2 + \varepsilon_1}{p_2 - \varepsilon_1} < 1 + \frac{\varepsilon}{p}$ (it is sufficient that $\varepsilon_1 < \frac{p_2\varepsilon}{2p+\varepsilon} $). Then, there exists $N_1>0$ such that 
$$
\left|n^{-2 r} \left( \frac{g_\Phi(n)}{g_\Phi(1/n)} - \frac{g_\Phi(n)}{r}\right)\right| < n^{-2r} \left(\frac{p_\Phi}{q_\Phi} - \frac{q_\Phi}{r}\right) <  \frac{\varepsilon_1}{2r}, \qquad n > N_1
$$
and
$$
p_2 - \frac{\varepsilon_1}{2} < g_\Phi(t) < p_2 + \frac{\varepsilon_1}{2}, \qquad t > N_1.
$$ 
Hence, if $l_n > 0$, then for $t>n>N_1$,
\begin{align*}
\dfrac{g_\Phi(t)}{1 + l_n \cdot(\Phi(t))^{-1}} &\le g_\Phi(t) < p + \varepsilon, \quad \text{and} \\
\dfrac{g_\Phi(t)}{1 + l_n \cdot(\Phi(t))^{-1}} 
&\ge \frac{g_\Phi(t)}{1 + l_n \cdot(\Phi(n))^{-1}} \\
&\ge \frac{r g_\Phi(t)}{p_2 + \varepsilon_1} \ge \frac{r (p_2 - \varepsilon_1)}{p_2 + \varepsilon_1} > q - \varepsilon.
\end{align*}
If instead $l_n < 0$, then
\begin{align*}
\dfrac{g_\Phi(t)}{1 + l_n \cdot(\Phi(t))^{-1}} &\ge g_\Phi(t) > q - \varepsilon, \quad \text{and} \\
\dfrac{g_\Phi(t)}{1 + l_n \cdot(\Phi(t))^{-1}} &\le \frac{g_\Phi(t)}{1 + l_n \cdot(\Phi(n))^{-1}} \\
&\le \frac{r g_\Phi(t)}{p_2 - \varepsilon_1} \le \frac{r (p_2 + \varepsilon_1)}{p_2 - \varepsilon_1} < p + \varepsilon.
\end{align*}
Choose $N_2>0$ such that
$$
p_1 - \varepsilon < g_\Phi(1/n) < p_1 + \varepsilon, \qquad n> N_2.
$$
If $k_n > 0$ then
$$
g_\Phi(1/n) \le \frac{r}{1 + k_n t^{-r}} < r, \qquad 1/n < t \le n,
$$
and if $k_n < 0$ then
$$
r< \frac{r}{1 + k_n t^{-r}} \le g_\Phi(1/n) \qquad 1/n < t \le n.
$$
Since $q-\varepsilon < r < p + \varepsilon$, we have
$$
q -\varepsilon < \frac{r}{1 + k_n t^{-r}} < p +\varepsilon, \qquad 1/n < t \le n,
$$
and, lastly, we note that 
$$
q - \varepsilon < g_\Phi(t) < p+\varepsilon, \qquad t\le 1/N_2.
$$
Let $n > \max\{N_1,N_2\}$. Then the arguments above now show that
\begin{align*}
p_{\Psi_n} &= \max \left\{\sup_{t\le 1/n} g_\Phi(t), \sup_{1/n < t \le n} \left(\dfrac{r}{1 + k_nt^{-r}}\right), \sup_{t>n} \left( \dfrac{g_\Phi(t)}{1 + l_n \cdot(\Phi(t))^{-1}}\right)\right\} \\
&< p + \varepsilon,
\end{align*}
and similarly, $q_{\Psi_n} > q-\varepsilon$. This completes the proof.
\end{proof}

With these results in mind, we can directly connect the limits of $\frac{t \Phi'(t)}{\Phi(t)}$ as $t$ approaches the origin or infinity to the equivalence class Lebesgue exponents as follows.

\begin{cor}\label{cor:p1p2equivpPhiqPhi}
Let $\Phi$ be a Young function with
\begin{align*}
\lim_{t\rightarrow0} \frac{t \Phi'(t)}{\Phi(t)} = p_0 
\quad&\text{and}\quad 
\lim_{t\rightarrow \infty} \frac{t \Phi'(t)}{\Phi(t)} = p_\infty. \\
\intertext{Then}
p_{[\Phi]} = \max \{p_0,p_\infty\} 
\quad&\text{and}\quad 
q_{[\Phi]} = \min\{p_0,p_\infty\}. 
\end{align*}
\end{cor}

\begin{proof}
Let $p = \max\{p_0,p_\infty\}$, $q=\min\{p_0,p_\infty\}$. 
It is evident that for any $\Psi\in [\Phi]$, 
$$
p_\Psi \ge p, \quad q_\Psi \le q,
$$ 
hence $p_{[\Phi]} \ge p$ and $q_{[\Phi]}\le q$. For any $\varepsilon>0$, let $\Psi_\varepsilon$ be chosen so that $p\le p_\Psi < p+\varepsilon $ and $q\ge q_\Psi > q-\varepsilon$ (as is possible by \cref{prop:p1p2pPhiqPhi}). Then
$$
p_{[\Phi]} = \inf\Sets{p_\Psi}{\Psi\in [\Phi]} \le \inf_{\varepsilon>0} p_{\Psi_\varepsilon} = p
$$
and
$$
q_{[\Phi]} = \sup\Sets{q_\Psi}{\Psi\in [\Phi]}  \ge \sup_{\varepsilon>0} q_{\Psi_\varepsilon} = q,
$$
giving the result.
\end{proof}

Using Lebesgue exponents for equivalence classes, we can improve \eqref{eq:qPhiincl} in the case that the limit of $\frac{t \Phi'(t)}{\Phi(t)}$ exists as $t$ approaches zero or infinity as follows.

\begin{prop}\label{prop:inclusionsEquiv}
Let $\Phi$ be a Young function and $[\Phi]$ its equivalence class. Suppose that 
$$
\lim_{t\rightarrow 0^+} \frac{t \Phi'(t)}{\Phi(t)} \quad \text{and} \quad
\lim_{t\rightarrow \infty} \frac{t \Phi'(t)}{\Phi(t)}
$$
exist. If $1 < q_{[\Phi]} \le p_{[\Phi]} < \infty$, then
$$
L^{p} \cap L^{q} \subseteq
L^\Phi \subseteq
L^{p} + L^{q}
$$
for all $p \in(p_{[\Phi]},\infty) $ and $q\in(1, q_{[\Phi]})$.
\end{prop}

\begin{proof}
For any $\Psi \in [\Phi]$,
$$
L^{p_\Psi} \cap L^{q_\Psi} \subseteq L^\Psi \subseteq L^{p_\Psi} + L^{q_\Psi}.
$$
Since $L^\Psi = L^\Phi$, the result now follows from \cref{lem:equivp1p2}, \cref{prop:p1p2pPhiqPhi}, and \cref{cor:p1p2equivpPhiqPhi}.
\end{proof}

Unfortunately, we can not guarantee that there exists a Young function $\Psi$ in the equivalence class of Young functions $[\Phi]$ such that $p_{\Psi} = p_{[\Phi]}$ and $q_\Psi = q_{[\Phi]}$, hence we can not guarantee that $L^{p_{[\Phi]}}\cap L^{q_{[\Phi]}}\subseteq L^\Psi$, nor that $L^\Phi \subseteq L^{p_{[\Phi]}} + L^{q_{[\Phi]}}$. In fact, the following provides a counter example.

\begin{example}
Let $\Phi(t) = t^2 \ln(2 + t)$. Then
$$
\lim_{t\rightarrow \infty} \frac{t \Phi'(t)}{\Phi(t)} = \lim_{t\rightarrow 0} \frac{t \Phi'(t)}{\Phi(t)} = 2
$$
hence $p_{[\Phi]} = q_{[\Phi]} = 2$ and so
$$
L^{2+\varepsilon} \cap L^{2-\varepsilon} \subseteq L^\Phi \subseteq L^{2+\varepsilon} + L^{2 - \varepsilon}
$$
for all $\varepsilon > 0$ by Proposition \ref{prop:inclusionsEquiv}. In fact, since $q_\Phi = 2$ we must have
$$
L^{2+\varepsilon}\cap L^{2} \subseteq L^\Phi \subseteq L^{2+\varepsilon} + L^2
$$
for every $\varepsilon>0$,
but $L^\Phi \neq L^2$, and there is no $\Psi\in[\Phi]$ with $p_\Psi = p_{[\Phi]} = 2$.
\end{example}

%%

% ------------------------------------------------------ %
\subsection{A different approach}\label{subsec:different}
% ------------------------------------------------------ %

So far, our attempts to analyze Orlicz spaces revolve around comparing the behaviors of Young functions (and in particular, $\Phi\in\Delta_2$) to that of $t^p$ for $1\le p < \infty$. One approach is to make this connection more literal and simply let the exponent $p$ depend on $t$. Somewhat naively, we write $\Phi(t) = \Phi(1) t^{r_\Phi(t)}$, for some function $\fdef{r_\Phi}{[0,\infty)}{[-\infty,\infty]}$. Then, by normalizing $\Phi(t)$ we obtain $\tilde
\Phi(t) = t^{r_\Phi(t)}$ for $\tilde \Phi(t) = \frac{\Phi(t)}{\Phi(1)}$, and therefore
$$
r_\Phi(t) = \frac{\ln \tilde\Phi(t)}{\ln t}.
$$
For analysis of indices related to this quantity, see \autocite{mali:1985}. Much in the same way as $g_\Phi(t) = \frac{t \Phi'(t)}{\Phi(t)}$,
this is a function which enables comparisons with Lebesgue spaces. If $\Phi(t) = t^p$, $1\le p < \infty$, then
$$
g_\Phi(t) = r_\Phi(t) = p.
$$
If $\Phi\in \Delta_2$, then by \cref{cor:normPsiprop}\ref{cor:propitem3}, we must have 
$$
q_\Phi \le r_\Phi(t) \le p_\Phi, \quad t>0.
$$
Therefore, if we let $a_\Phi = \inf_{t>0} r_\Phi(t)$ and $b_\Phi = \sup_{t>0} r_\Phi(t)$, then $q_\Phi \le a_\Phi \le b_\Phi \le p_\Phi$, and by definition of $r_\Phi$, we must have
\begin{align*}
\Phi(1) t^{a_\Phi} &\le \Phi(t) \le \Phi(1)t^{b_\Phi}, \qquad t\ge 1, \\
\Phi(1) t^{b_\Phi} &\le \Phi(t) \le \Phi(1)t^{a_\Phi}, \qquad t\le 1,
\end{align*}
which is (possibly) an improvement of \cref{cor:normPsiprop}\ref{cor:propitem3}.

Observe also that if we replace $\tilde \Phi$ with $c \tilde \Phi$ for any constant $c>0$ in the expression for $r_\Phi(t)$ above, then
$$
\frac{\ln c\tilde \Phi(t)}{\ln t} = \frac{c}{\ln t} + r_\Phi(t).
$$
On one hand, $\frac{c}{\ln t}$ is undefined for $t=1$. On the other hand, the term $c/\ln t$ vanishes as $t$ approaches the origin or infinity, hence it does not matter if $\Phi$ is normalized or not in the expression for $r_\Phi$ when considering only those limits.

In \cref{subsec:equivYfn}, we are not able to guarantee that the limit of $\frac{t\Phi'(t)}{\Phi(t)}$ as $t$ approaches zero or infinity is always well-defined. Looking instead at $r_\Phi(t)$, we can say the following.

\begin{prop}\label{prop:lnPhiexist}
Suppose that $\Phi\in\Delta_2$. Then the following is true.
\begin{enumerate}
\item \label{prop:lnitem1} The limits 
$$
r_0 = \lim_{t\rightarrow 0^+} r_\Phi(t) \quad \text{and}\quad 
r_\infty = \lim_{t\rightarrow\infty} r_\Phi(t)
$$
exist.

\vrum

\item \label{prop:lnitem2}The limits $r_0$ and $r_\infty$ fulfill
$$
q_\Phi \le \liminf_{t\rightarrow 0^+} g_\Phi(t) \le
r_0 \le \limsup_{t\rightarrow 0^+} g_\Phi(t) \le p_\Phi
$$
and
$$
q_\Phi \le \liminf_{t\rightarrow \infty} g_\Phi(t) \le
r_\infty \le \limsup_{t\rightarrow \infty} g_\Phi(t) \le p_\Phi,
$$
where $g_\Phi(t) =\frac{t\Phi'(t)}{\Phi(t)}$.
\vrum

\item \label{prop:lnitem3}If the limits
$$
\lim_{t\rightarrow 0^+} g_\Phi(t) 
\mathand
\lim_{t\rightarrow \infty} g_\Phi(t)
$$
also exist, then they equal $r_0$ and $r_\infty$, respectively.

\vrum

\item \label{prop:lnitem4}For every $\varepsilon>0$ there exists a $K_\varepsilon > 0$ such that
$$
t^{r_0 +\varepsilon} < \Phi(t) < t^{r_0 - \varepsilon}, \quad t < \frac{1}{K_\varepsilon}
$$
and
$$
t^{r_\infty -\varepsilon} < \Phi(t) < t^{r_\infty + \varepsilon}, \quad t > K_\varepsilon.
$$
\end{enumerate}
\end{prop}

\begin{proof}
It is evident that \ref{prop:lnitem2} implies \ref{prop:lnitem3}. We will show \ref{prop:lnitem1} in the case that $t$ tends to infinity, as the case that $t$ tends to zero follows by analogous argument. The parts \ref{prop:lnitem2} and \ref{prop:lnitem4} will be proven along the way.

Let
$$
p = \limsup_{t\rightarrow \infty} \frac{t \Phi'(t)}{\Phi(t)}
\quad \text{and} \quad q=\liminf_{t\rightarrow \infty} \frac{t \Phi'(t)}{\Phi(t)}.
$$
Since $\Phi\in\Delta_2$, we have $q\le p<\infty$. For any $\varepsilon>0$, there exists $L_\varepsilon >0$ such that
$$
q - \varepsilon<\frac{t \Phi'(t)}{\Phi(t)} < p +\varepsilon, \quad t > L_\varepsilon.
$$
Performing the same kinds of manipulations as in the proofs of \cref{prop:qpPhiIneq} and \cref{lem:equivp1p2} gives
$$
A_\varepsilon t^{q - \varepsilon} \le \Phi(t) \le B_\varepsilon t^{p+\varepsilon}, \quad t>L_\varepsilon
$$
where $A_\varepsilon =g(L_\varepsilon)$ for some increasing function $g$ and $B_\varepsilon = f(L_\varepsilon)$ for some decreasing function $f$.
This implies that
$$
\lim_{t \rightarrow \infty} \frac{\Phi(t)}{t^{p_1}} = 0 \mathand \lim_{t\rightarrow\infty} \frac{t \Phi(t)}{t^{p_2}} = \infty
$$
whenever $p_1 > p+\varepsilon$ and $p_2 < q-\varepsilon$, and since this holds for any $\varepsilon>0$, the statement is true for all $p_1 > p$ and $p_2 < q$. Since $\frac{\Phi(t)}{t^\alpha}$ is a decreasing, continuous function with respect to $\alpha$, there exists
$$
r = \inf \Sets{\alpha \in [q,p]}{\lim_{t\rightarrow \infty} \frac{\Phi(t)}{t^\alpha}= 0} = \sup \Sets{\alpha \in [q,p]}{\lim_{t\rightarrow \infty} \frac{\Phi(t)}{t^\alpha}= \infty}.
$$
Now, $\lim_{t\rightarrow\infty} \frac{\Phi(t)}{t^r}$ may be $0, \infty$, or may not exist at all, but for every $\varepsilon>0$,
$$
\lim_{t\rightarrow \infty} \frac{\Phi(t)}{t^{r+\varepsilon}} = 0 \quad \text{and} \quad \lim_{t\rightarrow\infty} \frac{\Phi(t)}{t^{r-\varepsilon}} = \infty.
$$
Hence, for every $j=1,2,3\dots$ there exists $M_{j,\varepsilon}>0$ and $N_{j,\varepsilon}$ such that
$$
\frac{\Phi(t)}{t^{r+\varepsilon}} < \frac{1}{j}, \quad t>M_{j,\varepsilon}
$$
and
$$
\frac{\Phi(t)}{t^{r-\varepsilon}}> j, \quad t>N_{j,\varepsilon}.
$$
By letting $K_\varepsilon = \max\{M_{1,\varepsilon}, N_{1,\varepsilon}\}$,
we obtain
\begin{equation}\label{proof:item4}
t^{r-\varepsilon}< \Phi(t) < t^{r+\varepsilon}, \quad t>K_\varepsilon.
\end{equation}
Taking logarithms, we now obtain
$$
r-\varepsilon < \frac{\ln \Phi(t)}{\ln t} < r + \varepsilon
$$
for all $t>\max\{K_\varepsilon,1\}$. Since $\varepsilon>0$ was arbitrary, we have proven \ref{prop:lnitem1} when $t$ approaches infinity by the discussion in the beginning of this section. Proceeding in the same way for the other case gives \ref{prop:lnitem1}. This also proves \ref{prop:lnitem2} by construction of $r$, and \eqref{proof:item4} shows \ref{prop:lnitem4}.
\end{proof}

As is shown above, for $\Phi,\Psi\in\Delta_2$, the limits of $\frac{\ln \Phi(t)}{\ln t}$ and $\frac{\ln \Psi(t)}{\ln t}$ as $t$ approaches zero or infinity are well-defined. If $\Phi$ and $\Psi$ are equivalent, then
$$
\frac{\ln \Phi(t)}{\ln t} \le \frac{\ln C\Psi(t)}{\ln t} = \frac{\ln C}{\ln t} + \frac{\ln \Psi(t)}{\ln t},
$$
and similarly
$$
\frac{\ln \Psi(t)}{\ln t} \le \frac{\ln C}{\ln t} + \frac{\ln \Phi(t)}{\ln t}.
$$
Taking limits as $t$ approaches zero or infinity will lead to the $\frac{\ln C}{\ln t}$ term vanishing, hence we obtain the following.

\begin{prop}\label{prop:lnequivPhi}
Suppose that $\Phi,\Psi\in\Delta_2$ are equivalent. Then
$$
\lim_{t\rightarrow 0^+} r_\Phi(t) = \lim_{t\rightarrow 0^+} r_\Psi(t) \quad \text{and}\quad
\lim_{t\rightarrow \infty} r_\Phi(t) = \lim_{t\rightarrow \infty} r_\Psi(t).
$$
\end{prop}

We can now prove the following analogue to \cref{prop:inclusionsEquiv}.

\begin{thm}\label{thm:lnPhiinclusions}
Suppose that $\Phi\in\Delta_2$. Let
\begin{align*}
r_0 &= \lim_{t\rightarrow 0^+} r_\Phi(t),
&r_\infty &= \lim_{t\rightarrow \infty} r_\Phi(t),
\end{align*}
and suppose that $r_0,r_\infty > 1$.
Then
$$
L^{p} \cap L^{q} \subseteq L^\Phi \subseteq L^{p} + L^{q}
$$
for every 
$$
p > \max\left\{r_0, r_\infty\right\}, 
\quad\text{and}\qquad
1<q < \min\left\{r_0,r_\infty \right\}.
$$
\end{thm}

\begin{proof}
Since $p_\Phi<\infty$ by assumption, $r_0,r_\infty < \infty$ by \cref{{prop:lnPhiexist}}\ref{prop:lnitem2}. Let
$$
p = \max\{r_0, r_\infty\}, \mathandd q = \min\{r_0,r_\infty\}.
$$
With \cref{prop:lnequivPhi} in mind, by the exact same construction as in the proof of \cref{prop:p1p2pPhiqPhi} (now with $r\in[q,p]$ with the $p$ and $q$ defined above), for every $q-1>\varepsilon>0$, there exists an equivalent Young function $\Psi_\varepsilon$ with $p\le p_{\Psi_\varepsilon} < p+\varepsilon$ and $q-\varepsilon < q_{\Psi_\varepsilon} \le q$.
Since
$$
L^{p_{\Psi_\varepsilon}} \cap L^{q_{\Psi_\varepsilon}} \subseteq L^{\Psi_\varepsilon} \subseteq L^{p_{\Psi_\varepsilon}} + L^{q_{\Psi_\varepsilon}}
$$
and $L^\Phi = L^{\Psi_\varepsilon}$, the result now follows.
\end{proof}

The benefit of \cref{thm:lnPhiinclusions} as opposed to \cref{prop:inclusionsEquiv} is that we do not need to assume that the limits $\lim_{t\rightarrow0^+} r_\Phi(t)$ and $\lim_{t\rightarrow \infty} r_\Phi(t)$ exist; they always exist by \cref{prop:lnPhiexist}\ref{prop:lnitem1}, assuming $\Phi$ fulfills the $\Delta_2$-condition.
% ----------------------------------------
%  Mixed Orlicz spaces
% ----------------------------------------
\section{Trace type Orlicz spaces}\label{sec:phi}
For mixed Orlicz spaces, we have similar inclusions to those of \cref{rem:qphi},  \cref{prop:inclusionsEquiv}, and \cref{thm:lnPhiinclusions}. To express these, we will need some (admittedly cumbersome) notation.

Let $\maclB_1$, $\maclB_2$ be Banach spaces. Then, formally, the space $(\maclB_1,\maclB_2)$ consists of all $\fdef{f}{\rdd}{\co}$ such that
$$
\nm{f_2}{\maclB_2} < \infty,
$$
where $f_2(y) = \nm{f(\cdot,y)}{\maclB_1}$.

With this notation, $L^{\Phi,\Psi} = (L^{\Phi},L^\Psi)$, and we observe that
$$
(L^{p_\Phi}\cap L^{q_\Phi}, L^{p_\Psi}\cap L^{q_\Psi}) \subseteq L^{\Phi,\Psi} \subseteq (L^{p_\Phi}+L^{q_\Phi}, L^{p_\Psi}+L^{q_\Psi}).
$$
If we let $\Phi = \Psi$, we obtain the spaces that we here refer to as trace type Orlicz spaces, and the inclusions above become
$$
(L^{p_\Phi}\cap L^{q_\Phi}, L^{p_\Phi}\cap L^{q_\Phi}) \subseteq L^{\Phi,\Phi} \subseteq (L^{p_\Phi}+L^{q_\Phi}, L^{p_\Phi}+L^{q_\Phi}).
$$
For the usual Lebesgue spaces, $(L^p,L^p) = L^p$ and the purpose of this section is to investigate the relationship between $L^{\Phi,\Phi}$ and $L^\Phi$. Generally, it is not the case that $(L^{p_\Phi}\cap L^{q_\Phi}, L^{p_\Phi}\cap L^{q_\Phi}) = L^{p_\Phi}\cap L^{q_\Phi}$ nor that $ (L^{p_\Phi}+L^{q_\Phi}, L^{p_\Phi}+L^{q_\Phi}) = L^{p_\Phi} + L^{q_\Phi}$, but
$$
(L^{p_\Phi}\cap L^{q_\Phi}, L^{p_\Phi}\cap L^{q_\Phi}) \subseteq L^{p_\Phi}\cap L^{q_\Phi} \subseteq L^\Phi \subseteq L^{p_\Phi} + L^{q_\Phi} \subseteq (L^{p_\Phi}+L^{q_\Phi}, L^{p_\Phi}+L^{q_\Phi}).
$$

To show that neither $L^{\Phi,\Phi} \subseteq L^\Phi$ nor $L^\Phi \subseteq L^{\Phi,\Phi}$ in general, we will provide examples which utilize the fact that if $p_\Phi \neq q_\Phi$, then
$$ 
L^{p_\Phi} + L^{q_\Phi} \subsetneq (L^{p_\Phi}+L^{q_\Phi}, L^{p_\Phi}+L^{q_\Phi}),
$$ and similarly
$$
(L^{p_\Phi}\cap L^{q_\Phi}, L^{p_\Phi}\cap L^{q_\Phi}) \subsetneq L^{p_\Psi}\cap L^{q_\Psi}.
$$

The first of these examples is as follows.

\begin{example}\label{ex:PhineqPhiPhi}
Let $d=2$ and let $\Phi(t) = t + t^2$. Then $L^\Phi(\rr 2) = L^1 \cap L^2$ and the norm $\nm{\cdot}{\Phi}$ is equivalent to
$$
\nm{f}{\Phi} = \nm{f}{L^1} + \nm{f}{L^2}.
$$
On the other hand, the norm for the space $L^{\Phi,\Phi}(\rr 2)$ then becomes
\begin{align*}
\nm{f}{\Phi,\Phi} 
&= \int\left(\int |f(x,y)| d x + \left(\int |f(x,y)|^2 d x\right)^{1/2} \right) d y \\
&+ 
\left(\int\left(\int |f(x,y)| d x + \left(\int |f(x,y)|^2 d x\right)^{1/2} \right)^2 d y\right)^{1/2}
\\
&= \nm{f}{L^{1,1}} + \nm{f}{L^{2,1}} + \nm{\big(\nm{f(\cdot,y)}{L^1} + \nm{f(\cdot,y)}{L^2}\big)}{L^2}.
\end{align*}
It is evident that $\nm{f}{\Phi} \le \nm{f}{\Phi,\Phi}$ for every $f\in L^{\Phi,\Phi}$.

Let $\{f_n\}_{n=2}^\infty$ be a sequence of functions with 
$$
f_n(x,y) = e^{-nx^2 + 2\sqrt{n-1}x y - y^2}, \qquad n=2,3,4,\dots.
$$ 
Then, since
$$
-nx^2 + 2\sqrt{n-1}x y - y^2 = -n\left(x - \frac{\sqrt{n-1}}{n} y \right)^2 - \left(1 - \frac{\sqrt{n-1}}{n} \right)^2 y^2,
$$
performing appropriate variable substitutions yields
$$
\nm{f_n}{\Phi} = \frac{\pi}{\sqrt{n - \sqrt{n-1}^2}} + \frac{\sqrt{\pi}}{\sqrt{2}\sqrt[4]{n -\sqrt{n-1}^2}} = \pi + \sqrt{\frac{\pi}{2}}, \qquad n=2,3,4,\dots.
$$
On the other hand, similar calculations yield
$$
\nm{f_n}{L^{2,1}} = \frac{\sqrt[4]{\pi}}{\sqrt[4]{2 n}} \frac{\sqrt{\pi}}{\sqrt{1 - \frac{n-1}{n}}} = C\sqrt[4]{n},
$$
where $C = \frac{\pi^{3 /4}}{2^{1/4}}$.

Hence, if we let
$$
g(x,y) = \sum_{n=2}\frac{1}{n^{5/4}}f_n(x,y),
$$
then
$$
\nm{g}{\Phi} \le \left(\pi + \sqrt \frac\pi 2\right) \sum_{n=3}^\infty \frac{1}{n^{5/4}} < \infty,
$$
but
$$
\nm{g}{\Phi,\Phi} \ge \nm{g}{L^{2,1}} = C \sum_{n=2}^\infty \frac{1}{n} = \infty.
$$
In other words, $g\in L^\Phi(\rr 2)$ but $g\notin L^{\Phi,\Phi}(\rr 2)$.
\end{example}

The example above shows that there exist Young functions $\Phi$ such that $L^\Phi \neq L^{\Phi,\Phi}$, and specifically $L^{\Phi,\Phi}\subsetneq L^\Phi$. One can perform similar arguments to show that more generally, if $\Phi(t) = t^q + t^p$ with $q<p<\infty$, then $L^{\Phi,\Phi}(\rdd) \subsetneq L^{\Phi}(\rd)$.

By duality arguments, we can take this example and transform it into an example where the opposite proper set inclusion holds.

\begin{example}
Pick $\Phi(t) = t^2 + t^3$ so that, by slightly altering the construction in \cref{ex:PhineqPhiPhi}, we have $L^{\Phi,\Phi} \subsetneq L^\Phi$. Since $p_\Phi = 3< \infty$, there exists a Young function $\Psi$ so that $L^\Psi = (L^\Phi)'$, the dual of $L^\Phi$ (cf. \autocite[Theorem 8.19]{adams:2003}). Moreover $L^{\Psi,\Psi} = (L^{\Phi,\Phi})'$ (see also \autocite[Sec. VII]{rao:1991}).
Since $(L^q\cap L^p)' = L^{q'} + L^{p'}$, where $q'$ and $p'$ fulfill
$$
\frac1q + \frac1{q'} = 1, \quad \frac1p + \frac1{p'} = 1,
$$
and $L^\Phi = L^2 \cap L^3$, we must have $L^\Psi = L^2 + L^{3/2}$. One choice of $\Psi$ is therefore given by
$$
\Psi(t) = 
\begin{cases}
t^2, & t\le 1 \\
2 t^{3/2} - 1, & t > 1
\end{cases}
$$
Since $L^{\Phi,\Phi} \subsetneq L^\Phi$, it therefore follows by duality that $L^\Psi \subsetneq L^{\Psi,\Psi}$.
\end{example}

% ---------------------------------------------------------- %
\subsection{The Young function modular}\label{subsec:modular}
% ---------------------------------------------------------- %

To compare the spaces $L^\Phi$ and $L^{\Phi,\Phi}$ in general, we will need to compare the norms $\nm{\cdot}{L^{\Phi}}$ and $\nm{\cdot}{L^{\Phi,\Phi}}$.
When possible, it is typically convenient to work with $\rho_\Phi(f) = \nm{\Phi(|f|)}{1}$ in place of $\nm{f}{\Phi}$ (in the same way that working with $\nm{|f|^p}{1} = \nm{f}{p}^p$ sometimes leads to simpler calculations than using $\nm{f}{p}$), despite the fact that $\rho_\Phi(f)$ is not generally a norm.  The following proposition relates the two quantities $\nm{f}{\Phi}$ and $\rho_\Phi(f)$ to each other.

\begin{prop}\label{prop:lPhirhoPhi}
Let $\Phi\in\Delta_2$. If $f$ is a  Lebesgue measurable function, then
$$
\nm{f}{\Phi} < \infty \iff \rho_\Phi(f)<\infty,
$$ 
and in this case, one of the three following statements is true:
\begin{enumerate}
\item \label{prop:lPhirhoPhiitem1}
$
1 < \nm{f}{\Phi} \le \rho_\Phi(f);
$

\vrum

\item \label{prop:lPhirhoPhiitem2}
$
\rho_\Phi(f) \le \nm{f}{\Phi}   < 1;
$

\vrum

\item \label{prop:lPhirhoPhiitem3}
$
\nm{f}{\Phi} = \rho_\Phi(f) = 1.
$
\end{enumerate}
\end{prop}

The first statement is a well-known result and is a direct consequence of, for instance, \autocite[p. 268, 8.9]{adams:2003}. 
Parts which can be inferred by the other statements, such as $\rho_\Phi(f) > 1$ if and only if $\nm{f}{\Phi} > 1$ are also commonly known (cf. \autocite[p. 45]{luxe:1955}). In an effort to be self-contained, we prove all parts of the Proposition here.

\begin{proof}
To abbreviate further the notations for this proof, let $\alpha = \nm{f}{\Phi}$ and $\beta = \rho_\Phi(f)$.

Since $\alpha = 0$ if and only if $f = 0$ if and only if (under the given assumptions on $\Phi$) $\beta = 0$, we perform the calculations of the proof under the assumption that $f\neq 0$, so that $\alpha,\beta > 0$. 

If $\alpha < \infty$ then $\nm{\tilde{f}}{\Phi} = 1$ for $\tilde{f} = f/\alpha$. This implies $\rho_\Phi(\tilde{f}) \le 1$, hence by \cref{prop:qpPhiIneq},
$$
\beta = \rho_\Phi(\alpha\tilde{f}) \le \alpha^p \rho_\Phi(\tilde{f}) \le \alpha^p < \infty,
$$
where $p=p_\Phi$ if $\alpha \ge 1$ and $p=q_\Phi$ otherwise.

If $\beta < \infty$ then, using the convexity of $\Phi$, for any $c\ge 1$,
$$
\rho_\Phi(f/c) \le c^{-1} \beta.
$$
Since $c^{-1}\beta \rightarrow 0$ as $c\rightarrow \infty$, there exists $c_0\ge 1$ such that  $c^{-1} \beta\le 1$ for all $c\ge c_0$. Then, by definition of $\alpha$, we must have $\alpha \le c_0$, hence $\alpha < \infty.$

To show that the three listed cases are the only possible outcomes, we begin by showing that $\alpha > 1$ if and only if $\beta > 1$, and that in this case $1 < \alpha \le \beta$. 

First suppose $\alpha > 1$. Then, since $\Phi$ is convex and $1/\alpha < 1$, $\Phi(|f|/\alpha) \le 1/\alpha \Phi(|f|)$. By definition of $\alpha$, 
$$
\rho_\Phi(f/\alpha) \le 1,
$$
but under the given assumptions on $\Phi$ (see \cref{rem:phidef}), we must have
$$
\rho_\Phi(f/\alpha) = 1,
$$
since otherwise there would exist $\lambda < \alpha$ with
$$
\rho_\Phi(f/\lambda) \le 1.
$$
Therefore, using the convexity of $\Phi$,
$$
1 = \rho_\Phi(f/\alpha) \le \frac{\beta}{\alpha},
$$
giving $1<\alpha \le \beta$.

On the other hand, if $\beta > 1$, then we must have $\alpha > 1$, since otherwise
$$
\rho_\Phi(f/\alpha) \ge \beta > 1
$$
($\Phi$ is increasing), but using the convexity of $\Phi$ again,
$$
\rho_\Phi(f/\beta)\le \frac{\beta}{\beta} = 1,
$$
hence $1 < \alpha \le \beta$ by definition of $\alpha$.

Now, we show that $\alpha < 1$ if and only if $\beta < 1$ and that in this case $\beta \le \alpha < 1$. (This will complete the proof, as the only other possible outcome will then be \ref{prop:lPhirhoPhiitem3}.)
Suppose $\alpha < 1$. Then $1/\alpha > 1$, and since $\Phi$ must be strictly increasing under the assumptions above,
$$
1 = \rho_\Phi(f/\alpha) > \beta.
$$
On the other hand, by \cref{prop:qpPhiIneq} yet again (noting that $q_\Phi \ge 1$),
$$
\rho_\Phi(f/\beta) \ge \beta^{1-q_\Phi} \ge 1 =\rho_\Phi(|f|/\alpha),
$$
so that $\beta\le \alpha < 1$.

Lastly, suppose $\beta < 1$. Then, repeating the last argument gives $\beta \le \alpha$ in the exact same way. Since $\beta < 1$, we are not in case \ref{prop:lPhirhoPhiitem1} and therefore, at least, $\alpha \le 1$. But $\alpha = 1$ must imply $\beta = 1$ by \cref{rem:phidef} again, hence the proof is complete.
\end{proof}

Using \cref{prop:qpPhiIneq} and arguing similarly to the proof above, we can make the following observation, which more closely mirrors the fact that $\nm{f}{p}^p = \nm{|f|^p}{1}$ for the usual Lebesgue spaces, with $p < \infty$.

\begin{cor}\label{cor:l1lPhi}
Suppose that $\Phi\in\Delta_2$ and $f\in L^\Phi$. Then
$$
\nm{f}{\Phi}^{p_\Phi} \le \rho_\Phi(f) \le \nm{f}{\Phi}^{q_\Phi} < 1 
$$
if $\nm{f}{\Phi} < 1$, and
$$
1 < \nm{f}{\Phi}^{q_\Phi} \le \rho_\Phi(f) \le \nm{f}{\Phi}^{p_\Phi}
$$
if $\nm{f}{\Phi} > 1$.
\end{cor}

One might wonder if in fact $\Phi(\nm{f}{\Phi}) = \rho_\Phi(f)$, as is the case when $\Phi(t) = t^{p}$, $p\ge1$. There are many counterexamples (see \cref{rem:end1}), and the following is one of them.

\begin{example}\label{ex:l1neqlPhi}
Let $\Phi(t) = t^2 + t^3$ and $f(x) = (1+|x|^2)^{1/2}$. Then
$$
\rho_\Phi(f) = 1 + \frac{\pi}{2} \approx 2.571
$$
but
$$
\nm{f}{\Phi} = \sqrt[3]{\frac12 + \sqrt{\frac14 - \frac{\pi^3}{6^3}}} + \sqrt[3]{\frac12 - \sqrt{\frac14 - \frac{\pi^3}{6^3}}} \approx 1.496,
$$
giving $\Phi(\nm{f}{\Phi}) \approx \Phi(1.496) \approx 5.586$.

On the other hand, \cref{cor:l1lPhi} tells us that
$\nm{f}{\Phi}^2 \le 1 + \frac{\pi}{2} \le \nm{f}{\Phi}^3$. Since $\nm{f}{\Phi}^2 \approx 2.238$ and $\nm{f}{\Phi}^3 \approx 3.348$, this provides a good estimate of $\rho_\Phi(f)$.
\end{example}

Of course, we are more likely to want to provide upper and lower bounds for $\nm{f}{\Phi}$ in terms of $\rho_\Phi(f)$, as this is generally a simpler quantity to compute. These bounds can easily be extracted from \cref{cor:l1lPhi} as follows.

\begin{cor}\label{cor:l1lPhi2}
Suppose that $\Phi\in\Delta_2$ and that $f$ is a Lebesgue measurable function with $\nm{f}{\Phi} < \infty$. Then
$$
\rho_\Phi(f)^{1/q_\Phi} \le \nm{f}{\Phi} \le \rho_\Phi(f)^{1/p_\Phi} \le 1 
$$
if $\rho_\Phi(f) < 1$, and
$$
1 \le \rho_\Phi(f)^{1/p_\Phi} \le \nm{f}{\Phi} \le \rho_\Phi(f)^{1/q_\Phi}
$$
if $\rho_\Phi(f) > 1$.
\end{cor}

Applying \cref{cor:l1lPhi2} to \cref{ex:l1neqlPhi}, we can see that
$$
\left(1+\frac{\pi}{2} \right)^{1/3} \le \nm{f}{\Phi} \le \left(1+\frac{\pi}{2} \right)^{1/2},
$$
and since $(1+\pi/2)^{1/3}\approx 1.370$, $(1+\pi/2)^{1/2} \approx 1.603$, we have arrived at a rough approximation for $\nm{f}{\Phi}$.

In the following proposition, we highlight the difference between $L^\Phi$ and $L^{\Phi,\Phi}$ in terms of the Young function modular $\rho_\Phi$.

\begin{prop}\label{prop:GfHf}
Let $\Phi\in\Delta_2$, let $\defrd$ be a function, and let
$$
G_f(y) = \Phi\big(\nm{f(\cdot,y)}{L^\Phi(\rd)}\big) 
\quad \text{and}\quad 
H_f(y) = \nm{\Phi(f(\cdot,y))}{L^1(\rd)}.
$$
Then
$$
\nm{f}{L^{\Phi,\Phi}(\rdd)} < \infty \iff \nm{G_f}{L^1(\rd)} < \infty,
$$
and
$$
\nm{f}{L^\Phi(\rdd)} < \infty \iff \nm{H_f}{L^1(\rd)} < \infty.
$$
\end{prop}

\begin{proof}
By \cref{prop:lPhirhoPhi},
$$
\nm{f}{\Phi,\Phi} = \Nm{\nm{f(\cdot,y)}{\Phi}}{\Phi} < \infty \iff \rho_\Phi\big( \nm{f(\cdot,y)}{\Phi}\big) = \nm{G_f}{L^1} < \infty.
$$
Similarly,
$$
\nm{f}{L^{\Phi}} < \infty \iff \rho_\Phi(f) = \nm{H_f}{L^1} < \infty.
$$
\end{proof}

% ---------------------------------------------------------------------------- %
\subsection{Sub- and supermultiplicative Young functions}\label{subsec:submult}
% ---------------------------------------------------------------------------- %

It is natural to consider sub- and supermultiplicative Young functions \autocite[ pp. 28]{rao:1991}, as such properties greatly simplify calculations with $\nm{\cdot}{\Phi}$. With respect to such considerations, we have the following.

\begin{prop}\label{lem:submult}
Let $\Phi\in\Delta_2$.
\begin{enumerate}
\item \label{lem:submultitem1}There exists a constant $C>0$ such that $\Phi(ab) \le C \Phi(a)\Phi(b)$ for all $a,b\ge0$ if and only if
$$
\rho_\Phi(f) \le C \Phi\left(\nm{f}{L^\Phi(\rd)}\right), \quad f\in L^\Phi(\rd). 
$$

\vrum

\item \label{lem:submultitem2} There exists a constant $C>0$ such that $\Phi(ab) \ge C \Phi(a)\Phi(b)$ for all $a,b\ge 0$, if and only if 
$$
\rho_\Phi(f) \ge C \Phi\left(\nm{f}{L^\Phi(\rd)}\right), \quad f\in L^\Phi(\rd). 
$$
\end{enumerate}
\end{prop}

\begin{proof}
Suppose that there exists a constant $C>0$ such that
$$
\Phi(ab) \le C \Phi(a) \Phi(b)
$$
for every $a,b\ge 0$. 
Since the desired result trivially holds whenever $\nm{f}{\Phi} = 0$, we assume that $\nm{f}{\Phi}\neq 0$.
Then, with $a = \nm{f}{\Phi}$ and $b = \frac{|f(x)|}{\nm{f}{\Phi}}$,
$$
\Phi(|f(x)|) \le C \Phi\left(\nm{f}{\Phi}\right) \Phi\left( \frac{|f(x)|}{\nm{f}{\Phi}}\right), \quad x\in\rd.
$$
By the definition of $\nm{\cdot}{\Phi}$,
$$
\int \Phi\left( \frac{|f(x)|}{\nm{f}{\Phi}}\right) d x \le 1.
$$
Hence
$$
\int \Phi(|f(x)|) \, d x \le C \Phi(\nm{f}{\Phi}) \cdot 1,
$$
meaning
$$
\rho_\Phi(f) \le C \Phi(\nm{f}{\Phi}), \quad f\in L^\Phi.
$$
This proves the first implication of \ref{lem:submultitem1}. To prove the same implication for \ref{lem:submultitem2}, follow the same steps, but use the fact (\cref{rem:phidef}) that 
$$
\int \Phi\left( \frac{|f|}{\nm{f}{\Phi}}\right) d x = 1.
$$
To prove the other implication of \ref{lem:submultitem1}, suppose that $\rho_\Phi(f) \le C \Phi\left(\nm{f}{\Phi}\right)$ for every $f\in L^\Phi$, but that $\Phi(ab) > C \Phi(a)\Phi(b)$ for some $a,b> 0$. Let
$$
g(x) = 
\begin{cases}
ab & x\in U,\\
0 & x\notin U,
\end{cases}
$$
where $U$ is a ball such that $\nm{g}{\Phi} = a$. (Since $\Phi\in \Delta_2$, this is fulfilled if $|U| = 1/\Phi(b)$.)

Then, for every $x\in U$,
$$
\Phi(|g(x)|) > C \Phi(\nm{g}{\Phi}) \Phi\left( \frac{|g(x)|}{\nm{g}{\Phi}}\right).
$$
Integrating both sides gives (since $g=0$ outside of $U$)
$$
\int \Phi(|g(x)|) \, dx > C \Phi(\nm{g}{\Phi}) \int \Phi \left( \frac{|g(x)|}{\nm{g}{\Phi}}\right) dx.
$$
Since $\Phi \in \Delta_2$ (\cref{rem:phidef} again), we have
$$
\int \Phi\left( \frac{|g(x)|}{\nm{g}{\Phi}}\right) d x = 1
$$
hence, we obtain
$$
\rho_\Phi(g) > C \Phi(\nm{g}{\Phi}),
$$
which is a contradiction. This proves \ref{lem:submultitem1}. The second implication of \ref{lem:submultitem2} follows analogously.
\end{proof}

To state the results in this section with more generality, we will consider the following two conditions.
Suppose that $\Phi, \Phi_1,\Phi_2$ are Young functions.
Then the conditions are as follows.

\begin{enumerate}[label={(C\arabic*)}]
\item \label{cond1} There exists a constant $C_1>0$ such that
\begin{equation*}
\Phi(C_1ab) \le \Phi_1(a) \Phi_2(b)
\end{equation*}
for all $a,b\ge 0$.

\vrum

\item \label{cond2} There exists a constant $C_2>0$ such that
\begin{equation*}
\Phi(C_2 a b) \ge\Phi_1(a) \Phi_2(b)
\end{equation*}
for all $a,b\ge 0$.
\end{enumerate}

\begin{rem}
If we set $\Phi = \Phi_1 = \Phi_2$ then \ref{cond1} and \ref{cond2} become
\begin{enumerate}[label={(C\arabic*$'$)}]
\item \label{cond1alt} There exists a constant $C_1>0$ such that
\begin{equation*}
\Phi(C_1ab) \le \Phi(a) \Phi(b)
\end{equation*}
for all $a,b\ge 0$.

\vrum

\item \label{cond2alt} There exists a constant $C_2>0$ such that
\begin{equation*}
\Phi(C_2 a b) \ge\Phi(a) \Phi(b)
\end{equation*}
for all $a,b\ge 0$.
\end{enumerate}

If $\Phi\in \Delta_2$, then \ref{cond1alt} is equivalent to submultiplicativity of $\Phi$ and  \ref{cond2alt} is equivalent to supermultiplicativity. 
\end{rem}

We note the following reformulation of Conditions \ref{cond1} and \ref{cond2}.

\begin{lem}\label{lem:equivsubmult}
Suppose that $\Phi$, $\Phi_1$, and $\Phi_2$ are bijective Young functions.
\begin{enumerate}
\item \label{equivitem1} Condition \ref{cond1} holds if and only if 
$$
C_1 \Phi_1^{-1}(a)\Phi_2^{-1}(b) \le \Phi^{-1}(ab)   
$$
for all $a,b\ge 0$.

\vrum

\item \label{equivitem2} Condition \ref{cond2} holds if and only if 
$$
\Phi^{-1}(ab)  \le  C_2\Phi_1^{-1}(a)\Phi_2^{-1}(b)   
$$
for all $a,b\ge 0$.
\end{enumerate}
\end{lem}

\begin{proof}
Suppose that \ref{cond1} holds. For any $a,b\ge 0$ let $x = \Phi_1^{-1}(a)$ and $y=\Phi_2^{-1}(b)$. Since
$$
\Phi(C_1 x y) \le \Phi_1(x) \Phi_2(y),
$$
we thus have
$$
\Phi(C_1 \Phi_1^{-1}(a) \Phi_2^{-1}(b)) \le a b.
$$
Since $\Phi$ is bijective, this inequality implies that
$$
C_1 \Phi_1^{-1}(a) \Phi_2^{-1}(b) \le \Phi^{-1}( a b),
$$
which shows one implication of \ref{equivitem1}.

Now suppose that 
$$
C_1 \Phi_1^{-1}(a) \Phi_2^{-1}(b) \le \Phi^{-1}( a b).
$$
for all $a,b\ge 0$. For any $a,b\ge 0$, let $x = \Phi_1(a)$ and $y = \Phi_2(b)$. Then, in particular,
$$
C_1 \Phi_1^{-1}(x)\Phi_2^{-1}(y) \le \Phi^{-1}(xy).
$$
Since $\Phi_1$ and $\Phi_2$ are bijective, we have
$$
C_1 a b \le \Phi^{-1}(\Phi_1(a) \Phi_2(b)).
$$
Since $\Phi$ is bijective,
$$
\Phi(C_1 a b) \le \Phi_1(a) \Phi_2(b),
$$
which gives the other implication of \ref{equivitem1}. The proof of \ref{equivitem2} follows analogously. 
\end{proof}

The following lemma follows by straight-forward computations.

\begin{lem}\label{lem:indnorm}
Suppose that $\Phi$, $\Phi_1$, and $\Phi_2$ are bijective Young functions. Let $A,B$ be measurable subsets of $\rd$ with $|A|>0$ and $|B|>0$ and let $f_{A,B}$ be the indicator function for the set $A\times B\subseteq \rdd$. Then
$$
\nm{f_{A,B}}{L^\Phi(\rdd)} = \frac{1}{\Phi^{-1}\left(\frac{1}{|A|\cdot|B|}\right)},
$$
$$
\nm{f_{A,B}(\cdot,y)}{L^{\Phi_1}(\rd)} =
\begin{cases}
\frac{1}{\Phi_1^{-1}\left(\frac{1}{|A|}\right)} & y\in B \\
0 & y\notin B,
\end{cases}
$$
and
$$
\nm{f_{A,B}}{L^{\Phi_1,\Phi_2}(\rdd)} = \frac{1}{\Phi_1^{-1}\left( \frac{1}{|A|}\right)\Phi_2^{-1}\left( \frac{1}{|B|}\right)}.
$$
\end{lem}

With this lemma in mind, we obtain the following proposition, which we prove here in an effort to be self-contained.

\begin{prop}\label{thm:general}
Suppose that $\Phi,\Phi_1, \Phi_2\in\Delta_2$. Then the following holds.
\begin{enumerate}
\item \label{generalitem1} We have $L^{\Phi_1,\Phi_2} \subseteq L^\Phi$ if and only if \ref{cond1} holds.

\vrum

\item \label{generalitem2} We have $L^{\Phi} \subseteq L^{\Phi_1,\Phi_2}$ if and only if \ref{cond2} holds.

\vrum

\item \label{generalitem3}  We have $L^{\Phi} = L^{\Phi_1,\Phi_2}$ if and only if $\Phi = \Phi_1 = \Phi_2$ and $L^\Phi = L^p$ for some $p\ge 1$.
\end{enumerate}
\end{prop}
\begin{proof}
Results \ref{generalitem1} and \ref{generalitem2} can be found in \autocite[p. 295]{mali:2004}. Moreover, \ref{generalitem3} is due to \autocite[p. 288-289]{fino:1991}. Here, we prove \ref{generalitem1} using a similar method to that of the latter proof. 

If $L^{\Phi_1,\Phi_2}\subseteq L^\Phi$, then there exists a constant $C>0$ such that
$$
\nm{f}{L^{\Phi}} \le C \nm{f}{L^{\Phi_1,\Phi_2}}, \quad f\in L^\Phi.
$$
Let $f_{A,B}$ be as in \cref{lem:indnorm}. By that same lemma,
$$
 \frac{1}{\Phi^{-1}\left(\frac{1}{|A|\cdot|B|}\right)} \le \frac{C}{\Phi_1^{-1}\left( \frac{1}{|A|}\right)\Phi_2^{-1}\left( \frac{1}{|B|}\right)}.
$$
For any $a,b> 0$, choose $A, B\subseteq \rd$ such that $|A| = \frac{1}{a}$ and $|B| = \frac{1}{b}$. Then
$$
\Phi_1^{-1} (a) \Phi_2^{-1}(b) \le C \Phi^{-1}(ab).
$$
Since the inequality trivially holds whenever $a=0$ or $b=0$, this is equivalent to \ref{cond1}, by \cref{lem:equivsubmult}.

Suppose instead that \ref{cond1} holds. Then for any $x,y\in \rdd$, by letting $a = \frac{|f(x,y)|}{\nm{f(\cdot,y)}{L^{\Phi_1}}}$ and $b=\nm{f(\cdot,y)}{L^{\Phi_1}}$ we obtain
$$
\Phi(|f(x,y)|) \le C_1 \Phi_1\left( \frac{|f(x,y)|}{\nm{f(\cdot,y)}{L^{\Phi_1}}}\right) \Phi_2(\nm{f(\cdot,y)}{L^{\Phi_1}}), \quad x,y\in\rd.
$$
Integrating with respect to $x$ and using the definition of $\nm{f(\cdot,y)}{L^{\Phi_1}}$, we obtain
$$
\int \Phi(|f(x,y)|) \, d x \le C \nm{f(\cdot,y)}{L^{\Phi_1}}.
$$
Integrating with respect to $y$ therefore gives
\begin{equation}\label{eq:modineq}
\rho_\Phi(f) \le C \rho_{\Phi_2}(\nm{f(\cdot,y)}{L^{\Phi_1}}).
\end{equation}
By \cref{prop:GfHf}, $\rho_\Phi(f) < \infty$ if and only if $f\in L^\Phi$ (since $p_\Phi < \infty$), and by slightly generalizing the statement of that same proposition, $\rho_{\Phi_2}(\nm{f(\cdot,y)}{L^{\Phi_1}}) < \infty$ if and only if $f\in L^{\Phi_1,\Phi_2}$ (since $p_{\Phi_2} < \infty$). Hence \eqref{eq:modineq} shows that $f\in L^{\Phi_1,\Phi_2}$ implies $f\in L^\Phi$. This proves \ref{generalitem1}. To prove \ref{generalitem2}, use analogous arguments with reverse inequalities, and utilize \cref{rem:phidef}.
\end{proof}

Setting $\Phi = \Phi_1 = \Phi_2$, we now obtain the following.

\begin{cor}\label{cor:tracetype}
Suppose that $\Phi\in \Delta_2$. Then the following holds.
\begin{enumerate}
\item \label{traceitem1}$L^{\Phi,\Phi} \subseteq L^\Phi$ if and only if $\Phi$ is submultiplicative.

\vrum

\item \label{traceitem2} $L^\Phi \subseteq L^{\Phi,\Phi}$ if and only if $\Phi$ is supermultiplicative.

\vrum

\item \label{traceitem3} $L^\Phi = L^{\Phi,\Phi}$ if and only if $L^\Phi = L^p$ for some $p\ge 1$.
\end{enumerate}
\end{cor}

\begin{rem}
Comparing \cref{cor:tracetype}\ref{traceitem1} and \ref{traceitem2} to \ref{traceitem3} of the same corollary, we see that $\Phi\in\Delta_2$ is both sub- and supermultiplicative if and only if $\Phi$ is equivalent to $t^p$ for some $p\ge 1$.
\end{rem}

Combining \cref{cor:tracetype} with \cref{lem:submult}, we immediately obtain the following.

\begin{cor}\label{cor:submult}
Suppose that $\Phi\in\Delta_2$. Then the following holds:
\begin{enumerate}
\item $L^{\Phi,\Phi} \subseteq L^\Phi$ if and only if there exists a constant $C>0$ such that
$$
\rho_\Phi(f) \le C \Phi\left(\nm{f}{L^\Phi}\right), \quad f\in L^\Phi;
$$

\vrum

\item $L^\Phi \subseteq L^{\Phi,\Phi}$ if and only if there exists a constant $C>0$ such that
$$
\rho_\Phi(f) \ge C \Phi\left(\nm{f}{L^\Phi}\right), \quad f\in L^\Phi; 
$$

\vrum

\item $L^\Phi = L^{\Phi,\Phi}$ if and only if there exists a constant $C>0$ such that
$$
C^{-1} \Phi\left(\nm{f}{L^\Phi}\right) \le \rho_\Phi(f) \le C \Phi\left(\nm{f}{L^\Phi}\right), \quad f\in L^\Phi,
$$
if and only if $\Phi$ is equivalent to $t^p$ for some $p\ge 1$.
\end{enumerate}
\end{cor}

\begin{rem}\label{rem:end1}
Evidently, if $\Phi(t) = t^p$, $1\le p < \infty$, then
$$
\Phi(\nm{f}{L^\Phi}) =  \int \Phi(|f(x)|) \, d x.
$$
On the other hand, if $\Phi(\nm{f}{L^\Phi}) = \int \Phi(|f(x)|)\, dx$ then $\nm{f}{L^\Phi}<\infty$ if and only if $\rho_\Phi(f)<\infty$ hence $\Phi\in \Delta_2$ (cf. \autocite[Section 8.9]{adams:2003}). In particular, $\Phi$ is bijective. Hence, for any $a,b>0$, with
$$
f(x) = \begin{cases}
    ab & x\in A \\
    0 & x\notin A,
\end{cases}
$$
where $A\subseteq \rd$ fulfills $|A| = \frac{1}{\Phi(b)}$, arguing as in \cref{lem:indnorm} we obtain
$$
\Phi(a)= \frac{\Phi(ab)}{\Phi(b)},
$$
which implies that $\Phi(t) = t^p$ for $p=\Phi'(1)$.
\end{rem}

Of course, we can generalize \cref{cor:submult} analogously to \cref{thm:general}, by generalizing \cref{lem:submult}. We state this result without proof, as it follows by analogous arguments.

\begin{thm}
Suppose that $\Phi,\Phi_1,\Phi_2\in \Delta_2$. Then the following holds:
\item $L^{\Phi_1,\Phi_2} \subseteq L^\Phi$ if and only if there exists a constant $C>0$ such that
$$
\rho_\Phi(f) \le C \Phi_2\left(\nm{f}{L^{\Phi_1}}\right), \quad f\in L^\Phi;
$$

\vrum

\item $L^\Phi \subseteq L^{\Phi_1,\Phi_2}$ if and only if there exists a constant $C>0$ such that
$$
\rho_\Phi(f) \ge C \Phi_2\left(\nm{f}{L^{\Phi_1}}\right), \quad f\in L^\Phi; 
$$

\vrum

\item $L^\Phi = L^{\Phi_1,\Phi_2}$ if and only if there exists a constant $C>0$ such that
$$
C^{-1} \Phi_2\left(\nm{f}{L^{\Phi_1}}\right) \le \rho_\Phi(f) \le C \Phi_2\left(\nm{f}{L^{\Phi_1}}\right), \quad f\in L^\Phi,
$$
if and only if $\Phi=\Phi_1=\Phi_2$ is equivalent to $t^p$ for some $p\ge 1$.
\end{thm}

%\includepdf[pages=-]{articles/AP23GS1.pdf}
% Biblatex
\printbibliography
% Bibtex
% \bibliographystyle{plain}
% \bibliography{ref}
\end{document}